\documentclass[12pt,letterpaper]{amsart}
\usepackage{amscd,amsxtra}

\newtheorem{theorem}{Theorem}[section]
\newtheorem{corollary}[theorem]{Corollary}
\newtheorem{lemma}[theorem]{Lemma}
\newtheorem{proposition}[theorem]{Proposition}

\theoremstyle{definition}

\newtheorem{remark}[theorem]{Remark}
\newtheorem{example}[theorem]{Example}

\DeclareMathOperator{\Dm}{\mathsf{D^{--}}}
\DeclareMathOperator{\Dcoh}{\mathsf{D_{coh}}}
\DeclareMathOperator{\Dbcoh}{\mathsf{D^{b}_{coh}}}
\DeclareMathOperator{\Dpcoh}{\mathsf{D^{+}_{coh}}}
\DeclareMathOperator{\Dmcoh}{\mathsf{D^{--}_{coh}}}
\DeclareMathOperator{\Dqc}{\mathsf{D_{qc}}}
\DeclareMathOperator{\Dbqc}{\mathsf{D^{b}_{qc}}}
\DeclareMathOperator{\Dpqc}{\mathsf{D^{+}_{qc}}}
\DeclareMathOperator{\Dmqc}{\mathsf{D^{--}_{qc}}}
\DeclareMathOperator{\FM}{\mathsf{FM}}
\newcommand{\ARG}{\,\cdot\,}
\newcommand{\cx}{^{\bullet}}
\newcommand{\dcx}{^{\prime\bullet}}
\DeclareMathOperator{\Picfu}{\textsf{Pic}}
\DeclareMathOperator{\CPic}{\overline{\textsf{Pic}}}
\DeclareMathOperator{\Coh}{\mathsf{Coh}}

\DeclareMathOperator{\Hilb}{\mathsf{Hilb}}
\DeclareMathOperator{\depth}{depth}
\DeclareMathOperator{\Spec}{\mathsf{Spec}}
\newcommand{\dtens}{\stackrel{\boldsymbol{L}}{\otimes}} 
\newcommand{\lar}{\longrightarrow}

\newcommand{\kL}{\mathcal L}

\newcommand{\kP}{\mathcal P}
\begin{document}

\title
[Fourier-Mukai transform on genus one fibrations]
{On a relative Fourier-Mukai transform on genus one fibrations} 

\author{Igor Burban}
\address{%
Max Planck-Institut f\"ur Mathematik,
Vivatsgasse 7, 53111 Bonn, Germany
}
\email{burban@mpim-bonn.mpg.de}

\author{Bernd Kreu{\ss}ler}
\address{%
Mary Immaculate College, South Circular Road, Limerick, Ireland
}
\email{bernd.kreussler@mic.ul.ie}

\subjclass[2000]{18E30, 14H60, 14H20, 14J27, 14H10}



\begin{abstract}

We study relative Fourier-Mukai transforms on genus one fibrations with
section, allowing explicitly the total space of the fibration to be singular
and non-projective. 
Grothendieck duality is used to prove a skew-commutativity
relation between this equivalence of categories and certain duality functors. 
We use our results to explicitly construct examples of semi-stable sheaves on
degenerating families of elliptic curves. 
\end{abstract}

\maketitle

\section{Introduction}
Mukai \cite{Mukai} introduced functors of the form
$\boldsymbol{R}\pi_{2\ast}(\mathcal{P}\dtens \pi_{1}^{\ast}(\ARG))$ as an
efficient tool to study vector bundles on Abelian varieties. If such a functor
is an equivalence of derived categories, it is called a Fourier-Mukai
transform, whereas in general they are referred to as integral transforms. 
Originally, they have been applied to study moduli problems on
smooth projective varieties, not only on complex tori. 
This seems to be natural in the light of Orlov's theorem \cite{Orlov}, which
states that any auto-equivalence of the bounded derived category of coherent
sheaves on a smooth projective variety $X$ is a Fourier-Mukai transform. 

More recently, in higher dimensional birational geometry, the point of view
was adopted that bounded derived categories of coherent sheaves might provide
the framework which is needed to understand the minimal model programme. 
This is supported by results which show that, if two threefolds are related by
a flop, their derived categories of bounded complexes of coherent sheaves
are equivalent, see \cite{BridgelandFlops, Chen, Kawamata}.

The minimal model programme naturally leads to the study of singular projective
varieties. With the recent development in mind, this generates a demand for
the study of Fourier-Mukai transforms on singular varieties.
However, not only Orlov's representability theorem, but most of the results
concerning equivalences between derived categories of coherent sheaves are
established in the smooth case only, see \cite{Orlov}.
The smoothness assumption enters the proofs in an essential way by using 
that any coherent sheaf has a finite locally free resolution and a finite
quasi-coherent injective resolution. Consequently, all the standard functors
of Grothendieck are defined on and take values in the bounded derived category
of coherent sheaves. 
This is no longer true on a singular variety. To deal with difficulties like
these, we apply the machinery which was developed in \cite{RD}. 

Using derived categories and Fourier-Mukai transforms, Chen \cite{Chen}
studied flops on three-dimensional projective varieties with terminal
Gorenstein singularities. 
Using non-commutative algebra, similar results were proved in a more general
setting by Van den Bergh \cite{vandenBergh}.
Some general properties of integral transforms on
singular varieties can be found in Chen's paper, but the main difficulties with
singularities are circumvented by embedding such a variety into a smooth
four-dimensional variety. Another method to circumvent these difficulties was
used by Kawamata \cite{Kawamata}.

On elliptically fibred smooth projective varieties, relative Fourier-Mukai
transforms have been applied successfully to the study of moduli of vector
bundles, see \cite{Bridgeland1,BridgeMaciocia1,Yoshioka}. Caldararu 
\cite{Caldararu} studied relative Fourier-Mukai transforms on smooth elliptic
threefolds and is forced to use twisted sheaves, because he does not suppose
the existence of a section. 
Relative Fourier-Mukai transforms on elliptic fibrations are now established
as important tools in other areas as well, such as string theory, where
D-branes are studied \cite{ACHY, DonagiPantev} and also in the study of
Calogero-Moser systems \cite{BenZviNevins}.

In this paper, we consider elliptic fibrations $q:X\rightarrow S$ with
irreducible fibres and with a section. 
We do not need to suppose $X$ or $S$ to be either projective or smooth. 
Under these assumptions, we show in Theorem \ref{thm:equiv} that a particular
integral transform is an auto-equivalence of derived categories of bounded
(resp.~bounded above) complexes of coherent sheaves on $X$. 
This generalises \cite{BBHM}, Theorem 2.8, {\cite{Bridgeland1}, Theorem 5.3
and \cite{BridgeMaciocia1}, Theorem 1.2. 

Because we allow $X$ to be singular, the requirement that the fibres are
irreducible is not too restrictive, at least in the two dimensional
case. Namely, if $X\rightarrow S$ is a smooth elliptic surface with a section,
it follows from Kodaira's classification
of singular elliptic fibres that the components of a singular
fibre, which are not met by the section, always form a negative definite
configuration. In fact, such a configuration can be contracted to a rational
double point.

Motivated by possible applications to mirror symmetry and integrable systems,
we look at this situation from the point of view of the  study of the
degeneration of derived categories in families of smooth  elliptic curves. 
If the elliptic curve degenerates to a singular irreducible curve, the total
space may be singular. Therefore, it is important to allow the total space $X$
to be singular.
In the singular case, however, it is difficult to prove that a given integral
transform is an equivalence of categories, because the methods of Bridgeland
\cite{BridgelandEquiv} and Bondal, Orlov \cite{BondalOrlov95} do not apply.
We overcome such difficulties by using a completely different strategy which
allows us to use results from our earlier paper \cite{BurbanKreussler}, in
which the special case $S=\Spec(\boldsymbol{k})$ was studied.

The plan of this article is the following. 
We start Section \ref{sec:fmt} with recalling basic properties of the
compactified relative Jacobian, following Altman and Kleiman
\cite{AltmanKleimanI, AltmanKleimanII}. 
Instead of using a Poincar\'e sheaf on the fibred 
product of the fibration $q:X\rightarrow S$ with its dual fibration, we prefer
to work on $X\times_{S}X$ and give an explicit description of the sheaf
$\mathcal{P}$ which defines the integral transform. Of course, these two
approaches are equivalent. We give detailed proofs of some important
properties of $\mathcal{P}$, versions of which can also be found in
\cite{BBHM}.
After proving a compatibility property between relative Fourier-Mukai
transforms and direct image functors of closed embeddings, we show our main
result in this section, Theorem \ref{thm:equiv}, which states that the
relative Fourier-Mukai transform with kernel $\kP$ is an auto-equivalence of
the derived category of coherent sheaves on the total space $X$.

In Section \ref{sec:duality}  we derive a certain skew-commutativity of the
relative Fourier-Mukai transform $\FM_{\mathcal{P}}$ and the derived 
functor $\boldsymbol{R}\mathcal{H}om(\ARG,\kL)$, where $\kL$ is a line
bundle on $X$. 
If $X$ is Gorenstein, this functor is a dualising functor on $X$ and we
obtained a generalisation of Mukai's result \cite{Mukai}, (3.8) as well as a
more general form of our result  \cite{BurbanKreussler}, Theorem 6.11.

In the final Section \ref{sec:appl}, we give two examples of flat families of
coherent sheaves on a singular fibration of cubics. In both examples, the
general fibre of the family is a semi-stable vector bundle on a smooth elliptic
curve. The degeneration on the singular fibre has an interesting semi-stable
direct summand. In one example it is a torsion free but not locally free
sheaf, whereas in the other example it is a vector bundle which is not a twist
of an Atiyah bundle.

\smallskip
\noindent
\textbf{Acknowledgement}. Both authors would like to thank
U.~Bruzzo, D.~Her\-n\'andez Ruip\'erez, B.~Keller, M.~Lehn and P.~Schapira for
fruitful discussions. The first-named author would like to thank R.~Rouquier
and the Institut de Math\'emati\-ques de Chevaleret in Paris for an invitation
and constant support. 
The main work on this article has been done during a visit
of the second-named author at the Institut de Math\'ematiques de Chevaleret,
Paris, which was made possible through Research Seed Funding at Mary
Immaculate College. Both authors are grateful to the Mathematisches
Forschungsinstitut Oberwolfach, where they could carry out part of the work on
this article.

\smallskip
\noindent
\textbf{Notation.} 
We fix an algebraically closed field $\boldsymbol{k}$ of characteristic
zero. Throughout this paper we work in the category of 
$\boldsymbol{k}$-schemes. Unless otherwise stated, a point is always a closed
point and a fibre means a fibre over a closed point. If $y\in Y$ is a point,
we denote by $\boldsymbol{k}(y)$ the residue field of $y$ and consider it as a
sheaf with support at $y$.
A morphism of schemes $f:Y\rightarrow T$ is said to have pure dimension $n$,
if $\dim \mathcal{O}_{Y_{t},y} = n$ for all points $y\in Y$, where $t=f(y)$
and $Y_{t}$ denotes the fibre of $f$ over $t\in T$. 
If $Y$ and $T$ are locally of finite type over $\boldsymbol{k}$, any morphism
$f:Y\rightarrow T$ is locally of finite type. Such a morphism is called
Cohen-Macaulay (resp.~Gorenstein), if it is flat and all fibres are
Cohen-Macaulay (resp.~Gorenstein) schemes. 

If $X$ is a separated Noetherian scheme of finite dimension, we denote by
$\Dcoh(X)$ the derived category of the category of complexes of
$\mathcal{O}_{X}$-modules whose cohomology sheaves are coherent. By
$\Dpcoh(X)$, $\Dmcoh(X)$, resp.~$\Dbcoh(X)$ we denote the full subcategories
of $\Dcoh(X)$ which consist of those objects whose cohomology vanishes in
sufficiently negative degrees, resp.~sufficiently positive degrees,
resp.~negative and positive degrees.  
Similarly, the notation $\Dqc(X)$, $\Dbqc(X)$, $\Dpqc(X)$, $\Dmqc(X)$ refers
to quasi-coherent cohomology.

It is well-known that $\Dbcoh(X)$ is equivalent to the derived category of the
category of bounded complexes of coherent sheaves and similarly for
$\Dmcoh(X)$, see e.g.~\cite{SGA}.

\section{Relative Fourier-Mukai transforms}\label{sec:fmt}

The main result of this section is Theorem \ref{thm:equiv}, which shows that a
certain relative Fourier-Mukai transform is an equivalence of categories in
two ways:
$$
\FM_{\mathcal{P}}^{-}:\Dmcoh(X) \rightarrow \Dmcoh(X)
\quad\text{and}\quad
\FM_{\mathcal{P}}^{\mathsf{b}} : \Dbcoh(X)\rightarrow \Dbcoh(X).
$$
In Remark \ref{rem:unbounded}, we explain how the techniques developed in
\cite{Spaltenstein, Lipman} for unbounded complexes can be used to extend this
result to obtain equivalences 
$$
\FM_{\mathcal{P}}^{+}:\Dpcoh(X) \rightarrow \Dpcoh(X)
\quad\text{and}\quad
\FM_{\mathcal{P}} : \Dcoh(X)\rightarrow \Dcoh(X).
$$
Our main result will be shown under the following assumptions:
\begin{equation*}
  \text{($\star$)}
  \begin{cases}
    &S \text{ and } X \text{ are reduced, connected and separated schemes;}\\ 
    &S \text{ is of finite type over }\boldsymbol{k};\\
    &q:X\rightarrow S \text{ is a flat and projective
      $\boldsymbol{k}$-morphism, whose fibres}\\ 
    &\text{are integral Gorenstein curves of arithmetic genus one;}\\
    &\text{there exists a section } \sigma:S\rightarrow X \text{ of $q$
      which factors through}\\
    &\text{the open set of points in $X$, at which $q$ is smooth.}
  \end{cases}
\end{equation*}
If these assumptions are satisfied, the image of the section $\sigma$ is a
Cartier divisor $\Sigma\subset X$. 
These assumptions imply that $S,X$ and $X\times_{S}X$ are Noetherian and of
finite dimension. 
By $\Delta\subset X\times_{S}X$ we denote the diagonal and by
$\mathcal{I}_{\Delta}$ its ideal sheaf. The diagonal embedding is denoted
$\delta: X \rightarrow X\times_{S}X$.

The sheaf $\mathcal{P}$ on $X\times_{S}X$,
which is used in the definition of $\FM_{\mathcal{P}}$, is
\begin{equation}
  \label{eq:poincare}
  \mathcal{P} := \mathcal{I}_{\Delta} \otimes \pi_{1}^{\ast}
  \mathcal{O}_{X}(\Sigma) \otimes \pi_{2}^{\ast} \mathcal{O}_{X}(\Sigma).
\end{equation}
The sheaf $\mathcal{P}$ is flat over both factors, because the diagonal has
this property and $\mathcal{O}_{X}(\Sigma)$ is locally free.

We use the following notation. If $q:X\rightarrow S$ is a morphism of schemes
and $s\in S, x\in X$ are points, the
fibre of $q$ over $s$ is denoted by $X_{s}$ and its embedding into $X$ by
$j_{s}:X_{s} \subset X$. We denote the two projections $X\times_{S}X
\rightarrow X$ by $\pi_{1}, \pi_{2}$ and the two projections $X_{s}\times
X_{s} \rightarrow X_{s}$ by $p_{1}, p_{2}$. Furthermore, we abbreviate $\pi :=
q\circ \pi_{1} = q\circ\pi_{2}: X\times_{S}X \rightarrow S$. 

The fibres of both projections $\pi_{1}$ and $\pi_{2}$ over $x\in X$ are
isomorphic to $X_{q(x)}$, so that we have two Cartesian squares: 
$$\begin{CD}
  X\times_{S}X @<{\psi_{2,x}}<< X_{q(x)}     @>{\psi_{1,x}}>> X\times_{S}X\\
  @V{\pi_{2}}VV                 @VVV                         @VV{\pi_{1}}V\\
  X            @<{x}<< \Spec(\boldsymbol{k}) @>{x}>>          X
\end{CD}$$
In addition, we have $\pi_{1}\circ \psi_{2,x} = \pi_{2}\circ \psi_{1,x} =
j_{q(x)}$. The morphisms $\psi_{\nu,x}$ coincide with the
compositions $$\psi_{\nu,x}: X_{q(x)}\rightarrow X_{q(x)}\times X_{q(x)}
\rightarrow X\times_{S}X,$$ where the first map embeds $X_{q(x)}$ as the fibre
of $p_{\nu}$ over $x$ and the second morphism is just $j_{q(x)}\times
j_{q(x)}$. By $X_{s} \cong \Delta_{s}\subset X_{s} \times X_{s}$ we denote the
diagonal and by $\mathcal{I}_{\Delta_{s}}$ its ideal sheaf.

Under the assumptions ($\star$) we obtain for any $s\in S$:
\begin{equation}
  \label{eq:fibre}
  (j_{s}\times j_{s})^{\ast}\mathcal{P} \cong \mathcal{I}_{\Delta_{s}} \otimes
  p_{1}^{\ast}\mathcal{O}_{X_{s}}(\sigma(s)) \otimes
  p_{2}^{\ast}\mathcal{O}_{X_{s}}(\sigma(s)).
\end{equation}
If $s=q(x)$, this implies 
\begin{equation}
  \label{eq:restr}
  \psi_{1,x}^{\ast}\mathcal{P} \cong \psi_{2,x}^{\ast}\mathcal{P} \cong
  \mathcal{I}_{x} \otimes \mathcal{O}_{X_{s}}(\sigma(s)),
\end{equation}
where $\mathcal{I}_{x}$ denotes the ideal sheaf of the point $x$ in
$X_{s}=X_{q(x)}$. Because $\sigma(s)$ is smooth in its fibre $X_{s}$, as a
special case we obtain 
\begin{equation}
  \label{eq:specrestr}
  \psi_{1,\sigma(s)}^{\ast}\mathcal{P} \cong
  \psi_{2,\sigma(s)}^{\ast}\mathcal{P} \cong 
  \mathcal{O}_{X_{s}}.
\end{equation}
The following result is an immediate consequence of
\cite[Theorem 1.9]{AltmanKleimanI}.

\begin{lemma}\label{basechange}
  Let $T$ be a scheme, which is locally of finite type over $\boldsymbol{k}$,
  $f:Y\rightarrow T$ a projective and flat morphism of schemes and
  $\mathcal{P}$ a $T$-flat coherent sheaf on $Y$. Suppose 
  $\mathcal{E}xt_{Y_{t}}^{i}(\mathcal{P}(t), \mathcal{O}_{Y_{t}})=0$ for all
  points $t\in T$ and all $i>0$, where $Y_{t}\subset Y$ is the fibre of
  $f$ over $t$ and $\mathcal{P}(t)$ denotes the restriction of $\mathcal{P}$ to
  $Y_{t}$. Then, the following holds:
  $$\mathcal{E}xt_{Y}^{i}(\mathcal{P},\mathcal{O})=0 \text{ for }
  i>0 \text{ and } \mathcal{P}^{\vee} \text{ is } T\text{-flat.}$$
  Furthermore, for any morphism $g:T'\rightarrow T$ there is an isomorphism
  $$(\boldsymbol{1}_{Y}\times g)^{\ast}(\mathcal{P}^{\vee}) \cong
  ((\boldsymbol{1}_{Y}\times g)^{\ast}\mathcal{P})^{\vee}$$ on $Y\times_{T}T'$.
\end{lemma}

\begin{corollary}\label{basechangecor}
  Let $f:Y\rightarrow T$ be a projective and flat morphism of schemes, with
  $T$ locally of finite type over $\boldsymbol{k}$. Suppose,
  $f$ is a Gorenstein morphism of pure dimension one and $\mathcal{P}$ is a
  $T$-flat coherent sheaf on $Y$ such that for any point $t\in T$ the
  restriction $\mathcal{P}(t)$ of $\mathcal{P}$ to the fibre $Y_{t}$ is torsion
  free. Then, the following holds: 
  $$\mathcal{E}xt_{Y}^{i}(\mathcal{P},\mathcal{O})=0 \text{ for }
  i>0 \text{ and } \mathcal{P}^{\vee} \text{ is } T\text{-flat.}$$
  Furthermore, for any morphism $g:T'\rightarrow T$ there is an isomorphism
  $$(\boldsymbol{1}_{Y}\times g)^{\ast}(\mathcal{P}^{\vee}) \cong
  ((\boldsymbol{1}_{Y}\times g)^{\ast}\mathcal{P})^{\vee}$$ on $Y\times_{T}T'$.
\end{corollary}

\begin{proof}
  Because, by assumption, $\mathcal{P}(t)$ is torsion free on the 
  curve $Y_{t}$, for any point $y\in Y_{t}$ we have
  $\depth(\mathcal{P}(t)_{y})\ge 1$. Because $Y_{t}$ has pure dimension one,
  this implies that $\mathcal{P}(t)_{y}$ is a maximal Cohen--Macaulay module
  over the local Gorenstein ring $\mathcal{O}_{Y_{t},y}$. Hence, by local
  duality \cite[Cor. 3.5.11]{BrunsHerzog}, we obtain
  $\mathcal{E}xt_{Y_{t}}^{i}(\mathcal{P}(t),\mathcal{O}_{Y_{t}})=0$ for all
  $i>0$. The claim follows now from Lemma \ref{basechange}.
\end{proof}

If the morphism $q:X\rightarrow S$ satisfies ($\star$), we can
apply the corollary to $f:Y\rightarrow T$ being either projection
$\pi_{1}, \pi_{2}:X\times_{S}X \rightarrow X$. 

In the course of the proof of the following proposition we use the
compactified Picard scheme as it was introduced by A.~Altman and S.~Kleiman in
\cite{AltmanKleimanI} and \cite{AltmanKleimanII}. Let us briefly recall some
notation and the results we are going to use.
Let $S$ be a scheme which is locally of finite type over $\boldsymbol{k}$ and
let $q:X\rightarrow S$ be a flat and projective morphism, whose fibres $X_{s}$
are integral curves.  
Without the assumption on the fibres, there are two functors 
$\Picfu_{X\mid S}^{-}$ and $\Picfu_{X\mid S}^{=}$ defined, which associate to 
any $S$-scheme the set of all equivalence classes of flat families of torsion
free rank one sheaves, which are, in addition, supposed to be Cohen-Macaulay
on the fibres in the second case. Two such families are equivalent, if they
differ by a twist with the pull-back of an invertible sheaf from the parameter
scheme. Because on an integral curve a torsion free sheaf is automatically
Cohen-Macaulay, these two functors coincide in our situation.  
We denote them by $\CPic_{X\mid S}$.

A pair $(P,\mathcal{F})$, which consists of an $S$-scheme $P$ and
a finitely presented sheaf $\mathcal{F}$ on $X\times_{S}P$, which is $P$-flat
and whose restrictions to fibres of the projection to $P$ are torsion free
sheaves of rank one, is said to represent the functor $\CPic_{X\mid S}$
if the following holds: for any $S$-scheme $T$ and any $T$-flat family
$\mathcal{G}$ of finitely presented torsion free sheaves of rank one on
$X\times_{S}T$ there exists a unique $S$-morphism $f:T\rightarrow P$,
such that there exists an invertible sheaf $\mathcal{A}$ on $T$ and an
isomorphism $(\boldsymbol{1}_{X}\times f)^{\ast}\mathcal{F} \cong
\mathcal{G}\otimes q_{T}^{\ast}\mathcal{A}$. By $q_{T}$ we denote here the
second projection $X\times_{S}T \rightarrow T$. 

In general, the functor $\CPic_{X\mid S}$ is not representable by an
$S$-scheme, because it is not a sheaf in the Zariski topology or in any finer
topology. Therefore, in 
general, this functor has to be sheafified. The sheafified functor in the
\'etale topology was studied in \cite{AltmanKleimanI}. However, under the
assumption of the existence of a section $\sigma:S\rightarrow X$ of $q$, which
factors through the smooth locus of $q$, it was shown in
\cite{AltmanKleimanII}, Theorem 3.4 (iii), that the functor $\CPic_{X\mid S}$
is a sheaf in the \'etale topology. Actually, they show it is a sheaf in the
finer fppf-topology. Hence, we can apply \cite{AltmanKleimanI}, Theorems 8.1
and 8.5, to conclude: $\CPic_{X\mid S}$ is representable by an $S$-scheme,
which is the disjoint union of infinitely many projective $S$-schemes 
$\overline{Pic}_{X\mid S}^{n}$, which represent the functors $\CPic_{X\mid
  S}^{n}$. The functor $\CPic_{X\mid S}^{n}$ is the open sub-functor of
$\CPic_{X\mid S}$ which parametrises sheaves of degree $n$ on the fibres.  
All these $S$-schemes $\overline{Pic}_{X\mid S}^{n}$ are isomorphic to each
other, because the tensor product with the invertible sheaf
$\pi_{1}^{\ast}\mathcal{O}_{X}(\Sigma)$ defines an isomorphism of functors
$\CPic_{X\mid S}^{n} \rightarrow \CPic_{X\mid S}^{n+1}$. 

Finally, assuming in addition that all fibres have arithmetic genus one, in
\cite{AltmanKleimanI}, Example 8.9 (iii), it was shown that the 
functor $\CPic_{X\mid S}^{-1}$ is represented by
$(X,\mathcal{I}_{\Delta})$. Hence, the functor $\CPic_{X\mid S}^{0}$ is
represented by any pair $(X,\mathcal{I}_{\Delta} \otimes
\pi_{1}^{\ast}\mathcal{O}_{X}(\Sigma) \otimes \pi_{2}^{\ast}\mathcal{A})$,
where $\mathcal{A}$ is an invertible sheaf on $X$. In particular, if the
assumptions ($\star$) are satisfied, the pair $(X,\mathcal{P})$ represents the
functor $\CPic_{X\mid S}^{0}$.

If $g:T\rightarrow S$ is a morphism of schemes, we obtain an isomorphism of
functors 
$\CPic_{X\mid S}^{0} \times_{S}T \cong \CPic_{X\times_{S}T\mid T}^{0}$. 
Hence, $(X\times_{S}T, (\boldsymbol{1}_{X}\times g)^{\ast}\mathcal{P})$
represents $\CPic_{X\times_{S}T\mid T}^{0}$. Here we used the isomorphism
$(X\times_{S}X)\times_{S}T \cong (X\times_{S}T) \times_{T} (X\times_{S}T)$.

\begin{proposition}\label{prop:i}
  Under the assumptions ($\star$) there exists a unique morphism
  $i:X\rightarrow X$ of $S$-schemes such that there exists an invertible sheaf
  $\mathcal{M}$ on $S$ and an isomorphism of sheaves on $X\times_{S}X$ 
  $$(\boldsymbol{1}_{X}\times i)^{\ast} \mathcal{P} \cong \mathcal{P}^{\vee} 
  \otimes \pi^{\ast}\mathcal{M}.$$
  Furthermore, $i$ is compatible with base change, which means the
  following. If we denote for any morphism $g:T\rightarrow S$ by $i_{g}:
  X\times_{S} T\rightarrow X\times_{S} T$ the map which is defined via
  universality by $\left((\boldsymbol{1}_{X\times_{S}X}\times g)^{\ast}
  \mathcal{P}\right)^{\vee}$ on $(X\times_{S}X)\times_{S}T$, then $i_{g} =
  i\times \boldsymbol{1}_{T}$. 
\end{proposition}

\begin{proof}
  As seen above, the pair $(X, \mathcal{P})$ represents the functor
  $\CPic_{X\mid S}^{0}$. Because $\mathcal{P}^{\vee}$ is
  $\pi_{2}$-flat and, by Corollary \ref{basechangecor},
  $\psi_{2,x}^{\ast}(\mathcal{P}^{\vee}) \cong
  (\psi_{2,x}^{\ast}(\mathcal{P}))^{\vee}$ is torsion free and of degree zero, 
  the universality of $(X, \mathcal{P})$ implies the existence of a unique
  morphism $i:X\rightarrow X$ of $S$-schemes which satisfies
  $(\boldsymbol{1}_{X}\times i)^{\ast} \mathcal{P} \cong \mathcal{P}^{\vee}
  \otimes \pi_{2}^{\ast}\mathcal{A}$ for some invertible sheaf $\mathcal{A}$ on
  $X$. Now, consider the Cartesian diagram, where $\sigma_{1}=(\sigma\circ q)
  \times \boldsymbol{1}_{X}$ 
  $$\begin{CD}
    X @>{\sigma_{1}}>> X\times_{S}X\\
    @V{q}VV              @V{\pi_{1}}VV\\
    S @>{\sigma}>>   X
  \end{CD}$$
  and restrict both sides of the isomorphism above to $X$ via
  $\sigma_{1}$. Because $(\boldsymbol{1}_{X}\times i)\circ
  \sigma_{1}=\sigma_{1}\circ i$ and $\pi_{2}\circ \sigma_{1} =
  \boldsymbol{1}_{X}$, this yields an isomorphism
  $i^{\ast}\sigma_{1}^{\ast}\mathcal{P} \cong
  \sigma_{1}^{\ast}\mathcal{P}^{\vee} \otimes \mathcal{A}$. Note that
  Corollary \ref{basechangecor} yields an isomorphism
  $\sigma_{1}^{\ast}(\mathcal{P}^{\vee}) \cong
  (\sigma_{1}^{\ast}\mathcal{P})^{\vee}$. But $\sigma_{1}\circ
  j_{s}=\psi_{1,\sigma(s)}$, hence $j_{s}^{\ast}\sigma_{1}^{\ast}\mathcal{P}
  \cong \psi_{1,\sigma(s)}^{\ast} \mathcal{P} \cong \mathcal{O}_{X_{s}}$ by
  (\ref{eq:specrestr}), and $\mathcal{A}$ must be trivial on the fibres of
  $q$. This implies that there exists an invertible sheaf $\mathcal{M}$ on
  $S$, such that $\mathcal{A} \cong q^{\ast}\mathcal{M}$.

  If $g:T\rightarrow S$ is an arbitrary morphism of schemes, from
  Corollary \ref{basechangecor} we obtain 
  $(\boldsymbol{1}_{X\times_{S}X}\times g)^{\ast}(\mathcal{P}^{\vee}) \cong
  ((\boldsymbol{1}_{X\times_{S}X}\times g)^{\ast}\mathcal{P})^{\vee}$.
  The definition of $i$, $\CPic_{X\times_{S}T\mid T}^{0} \cong \CPic_{X\mid
  S}^{0} \times_{S}T$ and the uniqueness of $i_{g}$ imply now the
  compatibility with base change.  
\end{proof}

\begin{remark}\label{restricti}
  Suppose, the morphism $g$ in the proposition is given by a point $s\in S$,
  considered as a morphism $s:\Spec(\boldsymbol{k})\rightarrow S$. Then, the
  induced morphism $i_{s}:X_{s}\rightarrow X_{s}$ coincides with the morphism
  $i$ used in \cite{BurbanKreussler}, if the reference point $p_{0}$ is chosen
  to be equal to $\sigma(s)$.
\end{remark}

\begin{lemma} Assuming ($\star$), we obtain
  $\mathcal{P}^{\vee\vee} \cong \mathcal{P}, i^{2} = \boldsymbol{1}_{X}$ and
  $i$ is an isomorphism.
\end{lemma}

\begin{proof}
  Here, we denote by $\mathcal{P}(x)$ the restriction of $\mathcal{P}$ to
  the fibre $X_{q(x)}$ of $\pi_{2}: X\times_{S}X\rightarrow X$ over $x\in
  X$. This is a torsion free sheaf on a Gorenstein curve, hence local duality
  applies to prove reflexivity of $\mathcal{P}(x)$. A local calculation shows
  that the canonical mappings
  $$\mathcal{P}(x)\rightarrow \mathcal{P}^{\vee\vee}(x)\rightarrow
  \mathcal{P}^{\vee}(x)^{\vee}$$ 
  and
  $$\mathcal{P}(x)\rightarrow \mathcal{P}(x)^{\vee\vee}\rightarrow
  \mathcal{P}^{\vee}(x)^{\vee}$$ 
  coincide. From Corollary \ref{basechangecor} we know the second arrow
  in both cases is an isomorphism. Now, we obtain that 
  the canonical isomorphism $\mathcal{P}(x)\rightarrow
  \mathcal{P}(x)^{\vee\vee}$ is isomorphic to the restriction of the canonical
  morphism of sheaves $\mathcal{P}\rightarrow\mathcal{P}^{\vee\vee}$. By
  Nakayama's lemma this morphism is surjective, whence we have an exact
  sequence
  \begin{equation}
    \label{seq:refl}
    0\rightarrow \mathcal{K} \rightarrow \mathcal{P} \rightarrow
  \mathcal{P}^{\vee\vee} \rightarrow 0.
  \end{equation}
  But, by Corollary \ref{basechangecor}, $\mathcal{P}^{\vee\vee}$ is
  $\pi_{2}$-flat, hence the restriction of (\ref{seq:refl}) to the fibre
  $X_{q(x)}$ is exact. Reflexivity of $\mathcal{P}(x)$ implies $\mathcal{K}(x)
  = 0$, and using Nakayama's lemma again, we obtain $\mathcal{K}=0$. Hence
  $\mathcal{P}$ is reflexive. The statement about $i^{2}$ is now an easy
  consequence of universality. 
\end{proof}

\begin{remark}
  It is interesting to observe that the sheaf $\mathcal{P}$ is Cohen-Macaulay,
  provided the assumptions ($\star$) are satisfied and $S$ itself is
  Cohen-Macaulay. This follows from \cite{BrunsHerzog}, Theorem 2.1.7, because
  $\mathcal{P}$ is Cohen-Macaulay on all fibres as seen in the proof of
  Corollary \ref{basechangecor}. By the same reason, $X$ is Cohen-Macaulay if
  $S$ has this property.
\end{remark}

Let $q:X\rightarrow S$ be a morphism which satisfies the conditions ($\star$). 
We are going to study several versions of integral transforms
$\FM_{\mathcal{Q}\cx}$, defined by the formula 
\begin{equation}
  \label{eq:FMT}
  \FM_{\mathcal{Q}\cx}(\mathcal{E}\cx)=
 \boldsymbol{R}\pi_{2\ast}(\mathcal{Q}\cx\dtens\pi_{1}^{\ast}(\mathcal{E}\cx)).
\end{equation}
If an integral transform $\FM_{\mathcal{Q}\cx}$ is an equivalence of
categories, it is called a Fourier-Mukai transform. The object $\mathcal{Q}\cx$
is called the \emph{kernel} of the integral transform $\FM_{\mathcal{Q}\cx}$.

\begin{remark}\label{rem:RD}
  Because $X\times_{S}X$ is Noetherian and of finite dimension, the functor
  $\pi_{2\ast}$ has finite cohomological dimension. Hence,
  $\boldsymbol{R}\pi_{2\ast}$ is way-out in both directions and \cite{RD}, I
  \S7, allows to define $\boldsymbol{R}\pi_{2\ast}$ 
  for unbounded complexes. Moreover, $\pi_{2\ast}$ is proper, hence we obtain
  $$\boldsymbol{R}\pi_{2\ast}:\Dcoh(X\times_{S}X) \rightarrow \Dcoh(X)$$
  with restrictions $\Dpcoh(X\times_{S}X) \rightarrow \Dpcoh(X),
  \Dmcoh(X\times_{S}X) \rightarrow \Dmcoh(X)$ and $\Dbcoh(X\times_{S}X)
  \rightarrow \Dbcoh(X)$. Similar remarks apply to the other projections which
  we use below. Moreover, from \cite{RD}, II.5.1, we obtain an isomorphism
  $$\boldsymbol{R}(\pi_{2\ast}\circ\pi_{23\ast}) \cong
  \boldsymbol{R}\pi_{2\ast} \circ \boldsymbol{R}\pi_{23\ast}$$
  of functors $\Dcoh(X\times_{S}X\times_{S}X) \rightarrow \Dcoh(X)$.

  To deal with unbounded complexes when deriving the tensor product is more
  difficult, because we normally don't have a finiteness condition as
  above. In \cite{RD}, II \S4, we find a definition of the derived tensor
  product $$\dtens:\Dmcoh(X) \times \Dmcoh(X) \rightarrow \Dmcoh(X).$$
  If $f:X\rightarrow Y$ is a morphism and $Y$ is locally Noetherian, there
  exists a derived functor $\boldsymbol{L}f^{\ast}:\Dmcoh(Y) \rightarrow
  \Dmcoh(X)$, \cite{RD}, II \S4. It satisfies
  $\boldsymbol{L}(f^{\ast}g^{\ast}) \cong \boldsymbol{L}f^{\ast}
  \boldsymbol{L}g^{\ast}$. Moreover, if $f$ has finite Tor-dimension,
  e.g.\/ if $f$ is flat, these are available for unbounded complexes with
  coherent cohomology.
  If $f:X\rightarrow Y$ is arbitrary and $\mathcal{F}\cx, \mathcal{G}\cx
  \in\Dm(Y)$, there is an isomorphism $\boldsymbol{L}f^{\ast}\mathcal{F}\cx
  \dtens \boldsymbol{L}f^{\ast}\mathcal{G}\cx 
  \cong \boldsymbol{L}f^{\ast}(\mathcal{F}\cx\dtens \mathcal{G}\cx)$,
  see \cite{RD}, II \S5.

  Finally, the projection formula
  $$\boldsymbol{R}f_{\ast}\mathcal{F}\cx \dtens \mathcal{G}\cx
  \cong
  \boldsymbol{R}f_{\ast}(\mathcal{F}\cx \dtens
  \boldsymbol{L}f^{\ast}\mathcal{G}\cx)$$
  is a very
  important tool in our calculations. In \cite{RD}, II \S5, this isomorphism
  is proved with $\mathcal{F}\cx\in\Dmcoh(X), \mathcal{G}\cx\in\Dmcoh(Y)$  and
  $f:X\rightarrow Y$ a quasi-compact morphism of Noetherian schemes of finite
  dimension. 

  With this preparation in mind, for any $\mathcal{Q}\cx\in
  \Dmcoh(X\times_{S}X)$ 
  it is now clear that the formula (\ref{eq:FMT}) defines a functor 
  $$\FM^{-}_{\mathcal{Q}\cx}:\Dmcoh(X) \rightarrow \Dmcoh(X).$$
  If $\mathcal{Q}\cx$ is isomorphic to a bounded complex of coherent
  $\pi_{1}$-flat sheaves with $\pi_{2}$-proper support, it was shown in
  \cite{Chen} that the restriction of $\FM_{\mathcal{Q}\cx}^{-}$ is a functor
  $$\FM_{\mathcal{Q}\cx}^{\mathsf{b}}:\Dbcoh(X)\rightarrow \Dbcoh(X).$$
  See Remark \ref{rem:unbounded} for a discussion how to extend the results of
  \cite{RD} and the definition of $\FM_{\mathcal{Q}\cx}$ to unbounded
  complexes.  
\end{remark}

\begin{example}\label{expl:FMT}
  Let $\mathcal{L}$ be a line bundle on $X$ and recall that
  $\delta:X\rightarrow X\times_{S}X$ denotes the diagonal embedding. In this
  case, $\delta_{\ast}\mathcal{L}$ is $\pi_{1}$-flat. Because
  $\boldsymbol{R}\pi_{2\ast} (\delta_{\ast}\mathcal{L}\dtens
  \pi_{1}^{\ast}\mathcal{E}\cx) 
  \cong \boldsymbol{R}\pi_{2\ast} (\delta_{\ast}(\mathcal{L}\otimes
  \boldsymbol{L}\delta^{\ast}\pi_{1}^{\ast}\mathcal{E}\cx))
  \cong \mathcal{L}\otimes \mathcal{E}\cx$,
  we obtain
  $$\FM_{\delta_{\ast}\mathcal{L}}^{-}(\mathcal{E}\cx) \cong
  \mathcal{L}\otimes \mathcal{E}\cx\quad \text{ and }\quad
  \FM_{\delta_{\ast}\mathcal{L}}^{\mathsf{b}}(\mathcal{E}\cx) \cong
  \mathcal{L}\otimes \mathcal{E}\cx.$$
  Here, we used $\pi_{2}\circ\delta = \boldsymbol{1}_{X} =
  \pi_{1}\circ\delta$ and the projection formula
  $$\delta_{\ast}\mathcal{L} \dtens \mathcal{G}\cx \cong
  \delta_{\ast}(\mathcal{L}\otimes\boldsymbol{L}\delta^{\ast}\mathcal{G}\cx)$$
  with $\mathcal{G}\cx=\pi_{1}^{\ast}\mathcal{E}\cx$.
  Note that $\delta_{\ast}$ and
  $\mathcal{L}\otimes(\ARG)$ are exact functors. In \cite{RD}, II.5.6, this
  formula is shown for $\mathcal{G}\cx\in\Dmcoh(X\times_{S}X)$. Therefore, we
  don't get from \cite{RD} the corresponding result for
  $\FM_{\delta_{\ast}\mathcal{L}}^{+}$ or
  $\FM_{\delta_{\ast}\mathcal{L}}$. See, however, Remark \ref{rem:unbounded}. 
\end{example}
Observe that $\FM^{-}_{\mathcal{Q}\cx}$ commutes with shift functors, that is
$$\FM_{\mathcal{Q}\cx}^{-}\circ[m] \cong [m]\circ \FM_{\mathcal{Q}\cx}^{-} 
\cong
\FM_{\mathcal{Q}\cx[m]}^{-}$$ for all $m\in\mathbb{Z}$. Moreover, the flat base
change theorem implies for any flat $S$-morphism $f:X\rightarrow X$ that there
exists an isomorphism of functors  $$f^{\ast}\circ \FM^{-}_{\mathcal{Q}\cx} 
\cong
\FM^{-}_{(\boldsymbol{1}_{X}\times f)^{\ast}\mathcal{Q}\cx}.$$
In the following lemmas we collected less obvious properties of
integral transforms, which are needed in the proof of our main result below.

\begin{lemma}\label{lem:fmtrsstr}
  Suppose $S$ is of finite type over $\boldsymbol{k}$, let $q: X\rightarrow S$
  be a flat and proper morphism of schemes and $\eta: Y\rightarrow X$ a closed
  subscheme. Suppose $\zeta: Z\rightarrow X$ is a closed subscheme such that
  the diagram 
  $$\begin{CD}
    Y\times_{S}Z @>{\eta\times\zeta}>> X\times_{S}X\\
    @V{pr_{1}}VV                       @V{\pi_{1}}VV\\
    Y            @>{\eta}>>            X
  \end{CD}$$
  is Cartesian, e.g.\/ $Z=q^{-1}(q(Y))$. Suppose $\mathcal{Q}\cx\in
  \Dmcoh(X\times_{S}X)$ and denote its derived restriction by
  $\mathcal{R}\cx=\boldsymbol{L}(\eta \times \zeta)^{\ast} \mathcal{Q}\cx$.
  Then, there is an isomorphism of functors $\Dmcoh(Y)\rightarrow\Dmcoh(X)$
  $$\FM_{\mathcal{Q}\cx}^{-} \circ \eta_{\ast} \cong \zeta_{\ast}\circ
  \FM_{\mathcal{R}\cx}^{-}.$$
\end{lemma}

\begin{proof}
  In the commutative diagram
  \begin{equation}\label{diag:FM}
  \begin{CD}
    Y @<{pr_{1}}<<  Y\times_{S}Z @>{pr_{2}}>>  Z\\
    @V{\eta}VV       @V{\eta\times\zeta}VV  @V{\zeta}VV\\
    X @<{\pi_{1}}<< X\times_{S}X @>{\pi_{2}}>> X
  \end{CD}
  \end{equation}
  the left square is Cartesian by assumption. Because $\pi_{1}$ and, hence,
  $pr_{1}$ are flat, from this Cartesian diagram we obtain an isomorphism of
  functors 
  \begin{equation}
    \label{eq:fmi}
    \pi_{1}^{\ast}           \circ  \eta_{\ast}   \cong
    (\eta\times\zeta)_{\ast} \circ  pr_{1}^{\ast}. 
  \end{equation}
  The four functors involved in the isomorphism (\ref{eq:fmi}) are exact,
  hence they coincide with their derived versions. Using that closed
  embeddings are proper, we obtain that (\ref{eq:fmi}) is an
  isomorphism of functors $\Dmcoh(Y) \rightarrow \Dmcoh(X\times_{S}X)$. 

  Because the schemes considered here are of finite dimension, and because
  $\eta\times\zeta$ is affine, which implies exactness of
  $(\eta\times\zeta)_{\ast}$, the projection formula gives a functorial
  isomorphism  
  \begin{equation}
    \label{eq:fmii}
     \mathcal{Q}\cx\dtens(\eta\times\zeta)_{\ast}\mathcal{F}\cx \cong
     (\eta\times\zeta)_{\ast}(\boldsymbol{L}(\eta\times\zeta)^{\ast}
     \mathcal{Q}\cx \dtens \mathcal{F}\cx)
  \end{equation}
  with $\mathcal{F}\cx\in\Dmcoh(Y\times_{S}Z)$ and $\mathcal{Q}\cx\in
  \Dmcoh(X\times_{S}X)$.  

  Finally, from commutativity of the right hand square in the diagram
  (\ref{diag:FM}) together with exactness of the functors $\zeta_{\ast}$ and
  $(\eta\times\zeta)_{\ast}$ we gain an isomorphism of
  functors $\Dmcoh(Y\times_{S}Z) \rightarrow \Dmcoh(X)$
  \begin{equation}
    \label{eq:fmiii}
    \boldsymbol{R}\pi_{2\ast} \circ (\eta\times\zeta)_{\ast} \cong
    \zeta_{\ast} \circ \boldsymbol{R}pr_{2\ast}.
  \end{equation}
  Putting these observations together, we obtain the desired
  isomorphism, which is functorial in $\mathcal{E}\cx\in \Dmcoh(Y)$.
  \begin{align*}
    \FM_{\mathcal{Q}}^{-}(\eta_{\ast} \mathcal{E}\cx)
      &\cong
      \boldsymbol{R}\pi{_{2\ast}}(\mathcal{Q}\cx
      \dtens\pi_{1}^{\ast}\eta_{\ast}\mathcal{E}\cx)
      \\
      &\cong 
      \boldsymbol{R}\pi{_{2\ast}}(\mathcal{Q}\cx \dtens
      (\eta\times\zeta)_{\ast}pr_{1}^{\ast} \mathcal{E}\cx)
      &\text{ using (\ref{eq:fmi})}
      \\
      &\cong
      \boldsymbol{R}\pi{_{2\ast}}((\eta\times\zeta)_{\ast}(
      \boldsymbol{L}(\eta\times\zeta)^{\ast}\mathcal{Q}\cx
      \dtens pr_{1}^{\ast}\mathcal{E}\cx))
      &\text{ using (\ref{eq:fmii})}
      \\
      &\cong
      \zeta_{\ast}\boldsymbol{R}pr_{2\ast}(\mathcal{R}\cx
      \dtens pr_{1}^{\ast}\mathcal{E}\cx)
      &\text{ using (\ref{eq:fmiii})}
      \\
      &\cong \zeta_{\ast}(\FM_{\mathcal{R}\cx}^{-}(\mathcal{E}\cx)). 
  \end{align*}     
\end{proof}

As a corollary, we obtain a useful and probably well-known result, a version
of which can be found in \cite{Chen}, Lemma 6.1.

\begin{corollary}\label{cor:fmtrestr}
  Suppose $S$ is of finite type over $\boldsymbol{k}$, $s\in S$ and
  $q:X\rightarrow S$ is flat and projective. Let $\mathcal{Q}\cx\in
  \Dmcoh(X\times_{S}X)$ and denote by 
  $\mathcal{Q}_{s}\cx :=\boldsymbol{L}(j_{s}\times j_{s})^{\ast}\mathcal{Q}\cx$
  its derived restriction. Then, there is an isomorphism of functors
  $\Dmcoh(X_{s})\rightarrow \Dmcoh(X)$
  $$j_{s\ast}\circ\FM_{\mathcal{Q}_{s}\cx}^{-} \cong
  \FM_{\mathcal{Q}\cx}^{-}\circ j_{s\ast}.$$ 
\end{corollary}

\begin{proof}
  This follows with $Y=Z=X_{s}$ and $\zeta=\eta=j_{s}$ from Lemma
  \ref{lem:fmtrsstr}. 
\end{proof}

\begin{lemma}\label{lem:tensor}
  Let $S$ be of finite type over $\boldsymbol{k}$, $X$ connected and reduced
  and $q:X\rightarrow S$ a flat, projective morphism.
  If $\mathcal{Q}\cx\in \Dbcoh(X\times_{S}X)$ is such 
  that $\FM_{\mathcal{Q}\cx}^{-}(\boldsymbol{k}(x)) \cong \boldsymbol{k}(x)$ 
  for all points $x\in X$, there exists a line bundle $\mathcal{L}$ on
  $X$ such that $\mathcal{Q}\cx\cong \delta_{\ast}\mathcal{L}$. This implies
  $\FM_{\mathcal{Q}\cx}^{-}(\mathcal{E}\cx)\cong
  \mathcal{E}\cx\otimes \mathcal{L}$, functorial in $\mathcal{E}\cx\in
  \Dmcoh(X)$. 
\end{lemma}

\begin{proof}
  By assumption, we know $\FM_{\mathcal{Q}\cx}^{-}(\boldsymbol{k}(x)) =
  \boldsymbol{R}\pi_{2\ast} (\mathcal{Q}\cx \dtens
  \pi_{1}^{\ast}\boldsymbol{k}(x))$ $\cong \boldsymbol{R}\pi_{2\ast}
  \psi_{1,x\ast} \boldsymbol{L}\psi_{1,x}^{\ast}\mathcal{Q}\cx \cong 
  j_{q(x)\ast} \boldsymbol{L}\psi_{1,x}^{\ast}\mathcal{Q}\cx$ is isomorphic
  to $\boldsymbol{k}(x)$. Hence,
  $\boldsymbol{L}\psi_{1,x}^{\ast}\mathcal{Q}\cx$ is a sheaf for all points
  $x\in X$. From \cite{BridgelandEquiv}, Lemma 4.3 we deduce now that
  $\mathcal{Q}\cx$ is a coherent sheaf $\mathcal{Q}$ on $X\times_{S}X$ which is
  $\pi_{1}$-flat. In particular, $\boldsymbol{L}\psi_{1,x}^{\ast}\mathcal{Q}\cx
  \cong \psi_{1,x}^{\ast}\mathcal{Q} \cong \boldsymbol{k}(x)$.

  Because $q$ is projective, there is a $q$-very ample line bundle
  $\mathcal{A}$ on $X$ and $\pi_{2}^{\ast}\mathcal{A}$ is $\pi_{1}$-ample. For
  large positive $m$ the canonical mapping
  \begin{equation}\label{eq:can}
    \pi_{1}^{\ast}\pi_{1\ast}(\mathcal{Q} \otimes
    \pi_{2}^{\ast}\mathcal{A}^{\otimes m}) \rightarrow \mathcal{Q} \otimes
    \pi_{2}^{\ast}\mathcal{A}^{\otimes m}
  \end{equation}
  is surjective and $R^{j}\pi_{1\ast}(\mathcal{Q} \otimes
  \pi_{2}^{\ast}\mathcal{A}^{\otimes m}) = 0$ for any $j>0$.

  From $H^{j}(\psi_{1,x}^{\ast}(\mathcal{Q}\otimes
  \pi_{2}^{\ast}\mathcal{A}^{\otimes 
  m})) \cong H^{j}(\psi_{1,x}^{\ast}\mathcal{Q} \otimes
  j_{q(x)}^{\ast}\mathcal{A}^{\otimes m}) \cong H^{j}(\boldsymbol{k}(x)
  \otimes j_{q(x)}^{\ast} \mathcal{A}^{\otimes m}) \cong
  H^{j}(\boldsymbol{k}(x))$, we obtain
  $$h^{j}(\psi_{1,x}^{\ast}(\mathcal{Q} \otimes
  \pi_{2}^{\ast}\mathcal{A}^{\otimes m})) =
  \begin{cases}
    1&  \text{ if } j=0,\\
    0&  \text{ if } j\ne 0.
  \end{cases}$$
  Hence, $\mathcal{B}:=\pi_{1\ast}(\mathcal{Q}\otimes
  \pi_{2}^{\ast}\mathcal{A}^{\otimes m})$
  is locally free of rank one on $X$. Therefore, from (\ref{eq:can}) we obtain
  a surjection 
  \begin{equation}
    \label{eq:quot}
    \mathcal{O}_{X\times _{S}X} \rightarrow \mathcal{Q}\otimes
    \pi_{1}^{\ast}\mathcal{B}^{\vee} \otimes
    \pi_{2}^{\ast}\mathcal{A}^{\otimes m}.
  \end{equation}
  Again, $\psi_{1,x}^{\ast}(\mathcal{Q}\otimes \pi_{1}^{\ast}\mathcal{B}^{\vee}
  \otimes \pi_{2}^{\ast}\mathcal{A}^{\otimes m}) \cong \boldsymbol{k}(x)$,
  which implies the existence of a unique $S$-morphism $\varphi:X \rightarrow
  \Hilb^{1}(X\mid S)$ such that (\ref{eq:quot}) is isomorphic to the pull-back
  $(\boldsymbol{1}_{X}\times
  \varphi)^{\ast}(\mathcal{O}_{X\times\Hilb^{1}(X\mid S)} \rightarrow
  \mathcal{U})$ of the universal quotient. But, by \cite{AltmanKleimanI}, Lemma
  8.7, $\Hilb^{1}(X\mid S)$ is
  represented by $(X,\mathcal{O}_{\Delta})$, hence $\varphi:X\rightarrow X$ 
  satisfies $(\varphi\times\boldsymbol{1}_{X})^{\ast}\mathcal{O}_{\Delta}
  \cong \mathcal{Q} \otimes \pi_{1}^{\ast}\mathcal{B}^{\vee} \otimes
  \pi_{2}^{\ast}\mathcal{A}^{\otimes m}$. Because (\ref{eq:quot}) restricted
  to the fibre of $\pi_{1}$ over $x$ is the canonical map
  $\mathcal{O}_{X_{q(x)}} \rightarrow \boldsymbol{k}(x)$, we must have
  $\varphi = \boldsymbol{1}_{X}$. Hence, $\mathcal{Q} \cong
  \mathcal{O}_{\Delta} \otimes \pi_{2}^{\ast}(\mathcal{B}\otimes
  \mathcal{A}^{\otimes -m})$. With $\mathcal{L}:= \mathcal{B}\otimes
  \mathcal{A}^{\otimes -m}$ the claim now follows. 
\end{proof}

Our main focus is on the integral transform whose kernel is the $\pi_{1}$-flat
sheaf $\mathcal{P}$ on $X\times_{S}X$, which was defined in
(\ref{eq:poincare}). In this case, we actually obtain a functor 
$$\FM_{\mathcal{P}}:\Dcoh(X) \rightarrow \Dcoh(X).$$
This follows because $\mathcal{P}\otimes \pi_{1}^{\ast}(\ARG)$ is an exact
functor and, as seen before, $\boldsymbol{R}\pi_{2\ast}$ is defined for
unbounded complexes. 
Moreover, the restriction to the full sub-categories $\Dbcoh(X),
\Dpcoh(X)$ and $\Dmcoh(X)$ define functors
\begin{align*}
  \FM_{\mathcal{P}}^{\mathsf{b}}&:\Dbcoh(X) \rightarrow \Dbcoh(X)\\
  \FM_{\mathcal{P}}^{+}&:\Dpcoh(X) \rightarrow \Dpcoh(X)\\
  \FM_{\mathcal{P}}^{-}&:\Dmcoh(X) \rightarrow \Dmcoh(X).
\end{align*}
In the definition of all of them, the derived tensor product is the ordinary
one, because $\mathcal{P}$ is $\pi_{1}$-flat.

\begin{theorem}\label{thm:equiv}
  If the assumptions ($\star$) hold,
  $$\FM_{\mathcal{P}}^{-}:\Dmcoh(X)\rightarrow\Dmcoh(X) \quad\text{ and }\quad
  \FM_{\mathcal{P}}^{\mathsf{b}}:\Dbcoh(X)\rightarrow\Dbcoh(X)$$ are
  equivalences. 
\end{theorem}

\begin{proof}
  A standard calculation, see \cite{Mukai, Chen, BondalOrlov95}, shows
  $$\FM_{\mathcal{P}}^{-}\circ\FM_{\mathcal{P}}^{-}\cong
  \FM_{\mathcal{Q}\dcx}^{-}
  \quad\text{ and }\quad
  \FM_{\mathcal{P}}^{\mathsf{b}} \circ \FM_{\mathcal{P}}^{\mathsf{b}} \cong
  \FM_{\mathcal{Q}\dcx}^{\mathsf{b}},$$  
  where
  $\pi_{ij}:X\times_{S}X\times_{S}X \rightarrow X\times_{S}X$ denote the
  projections and $\mathcal{Q}\dcx\cong
  \boldsymbol{R}\pi_{13\ast}(\pi_{12}^{\ast}\mathcal{P} \otimes
  \pi_{23}^{\ast}\mathcal{P})$ in both cases.
  All the necessary relations between derived functors, which are needed to
  calculate $\mathcal{Q}\dcx$, are contained in \cite{RD}, II \S5, provided we
  deal with  complexes which are bounded above. The main ingredient, which is
  not provided by \cite{RD} for complexes which are bounded below or unbounded
  is the projection formula, see Remark \ref{rem:unbounded}. 
  
  Observe that $\mathcal{Q}\dcx$ is an object of $\Dbcoh(X\times_{S}X)$ and
  that the tensor product is the usual one, because $\mathcal{P}$ is
  $\pi_{1}$-flat. With $$\mathcal{Q}\cx := (\boldsymbol{1}_{X}\times
  i)^{\ast}\mathcal{Q}\dcx[1] \in\Dbcoh(X\times_{S}X)$$ we obtain
  $$[1]\circ i^{\ast} \circ \FM_{\mathcal{P}}^{-} \circ \FM_{\mathcal{P}}^{-}
  \cong \FM_{\mathcal{Q}\cx}^{-}$$ and the same for $\FM^{\mathsf{b}}$.
  Our strategy is now to show $\FM_{\mathcal{Q}\cx}^{-}(\boldsymbol{k}(x))
  \cong \boldsymbol{k}(x)$ for all $x\in X$. This enables us to calculate
  $\mathcal{Q}\cx$ with the aid of Lemma \ref{lem:tensor}, hence gives us
  $\FM_{\mathcal{Q}\cx}^{\mathsf{b}}$ as well.   
  By Corollary \ref{cor:fmtrestr}, we have
  $\FM_{\mathcal{Q}\cx}^{-}\circ j_{s\ast}
  \cong
  [1] \circ i^{\ast} \circ \FM_{\mathcal{P}}^{-} \circ \FM_{\mathcal{P}}^{-}
  \circ j_{s\ast} \cong
  [1] \circ j_{s\ast} \circ  i_{s}^{\ast} \circ \FM_{\mathcal{P}_{s}}^{-} \circ
  \FM_{\mathcal{P}_{s}}^{-}$ for all $s\in S$.
  Now, because $\mathcal{P}$ is $\pi_{1}$-flat and $q$ is flat by assumption,
  the sheaf $\mathcal{P}$ is $\pi$-flat as well. Hence, the derived
  restriction $\mathcal{P}_{s}$ to the fibre of $\pi$ over $s\in S$ coincides
  with the usual restriction $(j_{s}\times j_{s})^{\ast}\mathcal{P}$. By
  (\ref{eq:fibre}), Remark \ref{restricti} and \cite{BurbanKreussler}, Remark
  2.17, we see that \cite{BurbanKreussler} applied to the curve $X_{s}$
  studies the functor $\mathbb{F}:= [1] \circ
  \FM_{\mathcal{P}_{s}}^{\mathsf{b}}$. This 
  implies that \cite{BurbanKreussler}, Theorem 2.18 gives an isomorphism
  $i_{s}^{\ast} \circ \FM_{\mathcal{P}_{s}}^{\mathsf{b}} \circ
  \FM_{\mathcal{P}_{s}}^{\mathsf{b}} \circ [1]
  \cong \boldsymbol{1}$. Therefore, we obtain
  \begin{equation}
    \label{eq:restrFMT}
    \FM_{\mathcal{Q}\cx}^{\mathsf{b}}\circ j_{s\ast} \cong j_{s\ast}.
  \end{equation}
  Hence, $\FM_{\mathcal{Q}\cx}^{-}(\boldsymbol{k}(x)) \cong
  \FM_{\mathcal{Q}\cx}^{\mathsf{b}}(\boldsymbol{k}(x)) \cong 
  \boldsymbol{k}(x)$ for all $x\in X$, and Lemma \ref{lem:tensor}
  implies the existence of a line bundle $\mathcal{L}$ on $X$ such that
  $\mathcal{Q}\cx\cong \delta_{\ast}\mathcal{L}$. Example \ref{expl:FMT}
  implies now $\FM_{\mathcal{Q}\cx}^{-}(\mathcal{E}\cx) \cong \mathcal{E}\cx
  \otimes \mathcal{L}$ and 
  $\FM_{\mathcal{Q}\cx}^{\mathsf{b}}(\mathcal{E}\cx) \cong \mathcal{E}\cx
  \otimes \mathcal{L}$. With  
  $\mathcal{E}\cx=j_{s\ast}\mathcal{O}_{X_{s}}$ and using (\ref{eq:restrFMT}),
  we see that $\mathcal{L}$ is trivial along the fibres of 
  $q:X\rightarrow S$, hence $\mathcal{L}\cong q^{\ast}\mathcal{N}$ with some
  line bundle $\mathcal{N}$ on $S$. 
  The functor $\mathbb{T}_{\mathcal{N}}$, which is the tensor product with
  $q^{\ast}\mathcal{N}$, is an equivalence on $\Dbcoh(X)$ and on $\Dmcoh(X)$
  and it commutes with $i^{\ast}$ and the shift functor $[-1]$. Therefore,
  $$\FM_{\mathcal{P}}^{-} \circ \FM_{\mathcal{P}}^{-} 
  \cong \mathbb{T}_{\mathcal{N}} \circ i^{\ast} \circ [-1],$$ which is an
  equivalence. The same statement is true for $\FM^{\mathsf{b}}$. This implies
  $\FM_{\mathcal{P}}^{-}$ and $\FM_{\mathcal{P}}^{\mathsf{b}}$ are
  equivalences as well.  
\end{proof}

\begin{remark}\label{rem:unbounded}
  The results, which are collected in Remark \ref{rem:RD} have been
  sufficient to prove that $\FM_{\mathcal{P}}^{-}$ and
  $\FM_{\mathcal{P}}^{\mathsf{b}}$ are equivalences. However, if we like to
  extend this result to $\FM_{\mathcal{P}}^{+}$ or $\FM_{\mathcal{P}}$, we
  need a definition of the derived tensor product, a formula for
  $\boldsymbol{L}(f^{\ast} g^{\ast})$ and the projection formula for unbounded
  complexes.
  With the exception of the projection formula, this was provided by
  Spaltenstein \cite{Spaltenstein}. A proof of the projection formula can be
  found in Lipman's notes \cite{Lipman}, Section 3.9. He deduces it from the
  classical version by using homological compatibility of
  $\boldsymbol{R}f_{\ast}$ and $\varinjlim$.

  The projection formula was used to show
  $\FM_{\mathcal{P}}\circ\FM_{\mathcal{P}} \cong \FM_{\mathcal{Q}\dcx}$ in the
  proof of Theorem \ref{thm:equiv}. It is also used in the proof of
  $\FM_{\delta_{\ast}\mathcal{L}}(\mathcal{E}\cx) \cong
  \mathcal{E}\cx\otimes\mathcal{L}$ and to show
  that $f^{\ast}\FM_{\mathcal{Q}\dcx}(\mathcal{E}\cx) \cong
  \FM_{(\boldsymbol{1}_{X}\times f)^{\ast}\mathcal{Q}\dcx}(\mathcal{E}\cx)$
  for unbounded complexes $\mathcal{E}\cx$ and flat morphisms $f$. 

  Because \cite{Lipman} is not yet available in final form, we confine
  ourselves to state the following extension of Theorem \ref{thm:equiv} in
  this remark only:
  $$\FM_{\mathcal{P}}^{+}:\Dpcoh(X)\rightarrow\Dpcoh(X) \quad\text{ and }\quad
  \FM_{\mathcal{P}}:\Dcoh(X)\rightarrow\Dcoh(X)$$ are equivalences of
  categories. The proof is the same as for the case of bounded above
  complexes, but is uses the results from \cite{Lipman}.
\end{remark}

\begin{remark}\label{rem:gap}
  In \cite{BBHM}, Theorem 2.8, it is claimed that
  $\FM_{\mathcal{P}}^{\mathsf{b}}$ is an auto-equivalence of the derived
  category $\Dbcoh(X)$, provided $X$ is a smooth K3 surface which is
  elliptically fibred over $S=\mathbb{P}^{1}$. We are not able to understand
  all details of that proof. To be more specific, in our notation, the authors
  argue that the square $\FM_{\mathcal{Q}\cx}^{\mathsf{b}} =
  \FM_{\mathcal{P}}^{\mathsf{b}} \circ \FM_{\mathcal{P}}^{\mathsf{b}}$ 
  of the relative Fourier-Mukai functor is isomorphic to the composition of
  the tensor product with 
  the shift of a certain line bundle $\mathcal{L}$ and an equivalence induced
  by the involution $i$ of the variety $X$.  The authors show that the complex
  $\mathcal{Q}\cx$ has cohomology only in degree one and construct a map $f:
  \mathcal{H}^1(\mathcal{Q}\cx) \rightarrow \zeta_*(\mathcal{L})$, where
  $\zeta = \boldsymbol{1}_{X} \times i : X \rightarrow X\times_{S}X$. For a
  sheaf $\mathcal{H}$ on $X\times_{S}X$ 
  and $z\in X\times_{S}X$ let us denote its fibre over $z$ by $\mathcal{H}(z)
  := \mathcal{H}\otimes \boldsymbol{k}(z)$. In \cite{BBHM} it is shown 
  that the induced map on fibres: $f(z): \mathcal{H}^1(\mathcal{Q}\cx)(z) \lar
  \zeta_*(\mathcal{L})(z)$ is an isomorphism for all $z\in X\times_S X$.  This
  would be enough to prove that $\mathcal{H}^1(\mathcal{Q}\cx)$ is isomorphic
  to $\zeta_*(\mathcal{L})$, if it was known that
  $\mathcal{H}^1(\mathcal{Q}\cx)$ is a sheaf of
  $\mathcal{O}_{\Gamma}$-modules, where $\Gamma=\zeta(X)$ is the graph of
  $i$. However, we do not see an easy argument to show this. 
\end{remark}

\section{Fourier-Mukai transform and Grothendieck duality}\label{sec:duality}

In this section we study the relationship between duality and the relative
Fourier-Mukai transform $\FM_{\mathcal{P}}$ which was studied in the previous
section. 
We shall use Grothendieck-Verdier duality to prove our result, which describes
the composition of the Fourier-Mukai transform with a functor of the type
$\boldsymbol{R}\mathcal{H}om(\ARG, \mathcal{L})$. A corollary of our result
can be formulated as a kind of twisted compatibility with a dualising functor,
but only if we assume the scheme $S$ to be Gorenstein. As a special case, we
recover \cite{BurbanKreussler}, Theorem 6.11. 

The derived functor of $\mathcal{H}om_{X}(\ARG,\ARG)$ is defined by replacing
the second argument with an injective resolution. The result are bi-functors
\begin{align*}
  \boldsymbol{R}\mathcal{H}om_{X}(\ARG,\ARG)&:
  \Dmcoh(X)\times\Dpcoh(X)\rightarrow \Dpcoh(X) \quad\text{ and }\\
  \boldsymbol{R}\mathcal{H}om_{X}(\ARG,\ARG)&:
  \Dcoh(X)\times\Dbcoh(X)_{\mathsf{fid}}\rightarrow \Dcoh(X),
\end{align*}
where $\Dbcoh(X)_{\mathsf{fid}}$ denotes the full subcategory of $\Dbcoh(X)$
consisting of complexes with finite injective dimension, see \cite{RD}, II.3.3.
If $\mathcal{F}\cx\in\Dmcoh(X), \mathcal{G}\cx\in\Dpcoh(X)$ and $\mathcal{L}$
is a locally free sheaf on $X$, by \cite{RD}, II.5.16, we have a functorial
isomorphism 
$$\boldsymbol{R}\mathcal{H}om_{X}(\mathcal{F}\cx,
\mathcal{G}\cx\otimes\mathcal{L}) 
\cong
\boldsymbol{R}\mathcal{H}om_{X}(\mathcal{F}\cx, \mathcal{G}\cx) 
\otimes\mathcal{L}.$$
Because we allow $X$ and $S$ to be singular, we cannot expect that the functor
$\boldsymbol{R}\mathcal{H}om_{X}(\ARG,\mathcal{O}_{X}) : \Dmcoh(X) \rightarrow
\Dpcoh(X)$ sends $\Dbcoh(X)$ to $\Dbcoh(X)$. 
This would be true, if we supposed $S$ to be Gorenstein, but not in general. 
This forces us to carefully check the availability in the case of unbounded
complexes of standard results like the base change theorem and the projection
formula.
Throughout this section, we need to use the projection formula only
in the form in which it was presented in \cite{RD}.

In addition to the statements collected in
Remark \ref{rem:RD}, which are established in \cite{RD}, we need
the projection formula 
$$\boldsymbol{R}f_{\ast}(\mathcal{F}\cx)\dtens\mathcal{G} \cong 
\boldsymbol{R}f_{\ast}(\mathcal{F}\cx\dtens f^{\ast}\mathcal{G})$$ 
in the following situation:
$f:X\rightarrow Y$ is a flat and proper morphism of Noetherian schemes of
finite dimension, $\mathcal{F}\cx\in\Dpcoh(X)$ and $\mathcal{G}$ is a locally
free sheaf on $Y$.  
In this case, $\ARG\otimes \mathcal{G}$ and $\ARG \otimes f^{\ast}\mathcal{G}$
are both exact. The expressions on both sides are defined by 
replacing $\mathcal{F}\cx$ with an injective resolution. The usual projection
formula for sheaves provides the morphism and with the way-out technique from
\cite{RD}, I.7.1 (ii), it follows that we obtained an isomorphism. 

In addition to these standard formulas, we shall need in our proof of Theorem
\ref{thm:alldual} an isomorphism, which is provided in the following lemma. 

\begin{lemma}\label{lem:dual}
  Let $T$ be a scheme which is of finite type over $\boldsymbol{k}$
  and let $f:Y\rightarrow T$ be a projective and flat morphism whose fibres
  are Gorenstein curves. Suppose 
  $\mathcal{P}\in\Coh(Y)$ is a $T$-flat sheaf, which is torsion free on
  fibres. Then there exists an isomorphism in $\Dpcoh(Y)$
  $$\boldsymbol{R}\mathcal{H}om_{Y}(\mathcal{P} \otimes f^{\ast}\mathcal{E}\cx,
  \mathcal{O}) \cong
  \mathcal{P}^{\vee} \otimes
  f^{\ast}\boldsymbol{R}\mathcal{H}om_{T}(\mathcal{E}\cx,\mathcal{O}),$$
  which is functorial in $\mathcal{E}\cx\in \Dmcoh(T)$.
\end{lemma}

\begin{proof}
  The assumptions allow us to apply Corollary \ref{basechangecor} to obtain
  flatness of $\mathcal{P}^{\vee}$ over $T$ and 
  $\mathcal{E}xt_{Y}^{i}(\mathcal{P},\mathcal{O})=0 \text{ for } i>0.$

  For any open subset $U\subset T$ we have a Cartesian diagram
  $$\begin{CD}
    Y_{U} @>{j_{U}}>> Y\\ @VV{f_{U}}V @VV{f}V\\ U @>{i_{U}}>> T.
  \end{CD}$$
  We denote $\mathcal{P}_{U}:=j_{U}^{\ast}(\mathcal{P})$ and observe
  $\mathcal{P}_{U}^{\vee} \cong j_{U}^{\ast}(\mathcal{P}^{\vee})$. For any open
  $U\subset T$ we study the two functors 
  \begin{align*}
    \mathbb{A}_{U} &:=
    \boldsymbol{R}\mathcal{H}om_{Y_{U}}(
    \mathcal{P}_{U} \otimes f_{U}^{\ast}(\ARG),
    \mathcal{O}_{Y_{U}}) \quad\text{ and }\\
    \mathbb{B}_{U} &:=
    \mathcal{P}_{U}^{\vee} \otimes f_{U}^{\ast}
    \boldsymbol{R}\mathcal{H}om_{U}(\ARG, \mathcal{O}_{U}).
  \end{align*}
  Because $\mathcal{P}$ is supposed to be $T$-flat, $\mathcal{P}_{U}$ is
  $U$-flat and the functor $\mathcal{P}_{U} \otimes f_{U}^{\ast}(\ARG)$ is
  exact. From \cite{RD}, I \S 7, it is, therefore, clear that both functors
  $\mathbb{A}_{U}$ and $\mathbb{B}_{U}$ are contravariant way-out right
  $\partial$-functors. Because the sheaf $\mathcal{E}xt^{i}(\mathcal{F},
  \mathcal{G})$ is coherent for any $i$ and arbitrary coherent sheaves
  $\mathcal{F}, \mathcal{G}$, it follows from \cite{RD}, II.7.3, that we have
  $$\mathbb{A}_{U}, \mathbb{B}_{U}: \Dmcoh(U) \rightarrow \Dpcoh(Y_{U}).$$
  Our main aim is to show $\mathbb{A}_{T} \cong \mathbb{B}_{T}$. We are going
  to apply the lemma on way-out functors \cite{RD}, Prop.\/ I.7.1.
  This requires that we have well-defined natural transformations
  $\eta_{U}:\mathbb{A}_{U} \rightarrow \mathbb{B}_{U}$  for any open $U\subset
  T$ and a compatibility between them.
  For any $\mathcal{E}\cx\in\Dmcoh(U)$ we define $$\eta_{U}(\mathcal{E}\cx):
  \mathbb{A}_{U}(\mathcal{E}\cx) \rightarrow \mathbb{B}_{U}(\mathcal{E}\cx)$$
  as follows.
  We choose a resolution by injective quasi-coherent sheaves $\mathcal{O}_{U}
  \rightarrow 
  \mathcal{I}_{U}\cx$. Because $f_{U}$ is flat, $\mathcal{O}_{Y_{U}}
  \cong f_{U}^{\ast}\mathcal{O}_{U} \rightarrow
  f_{u}^{\ast}\mathcal{I}_{U}\cx$ is still a quasi-isomorphism.
  Let now $\mathcal{O}_{Y_{U}} \rightarrow \mathcal{J}_{U}\cx$ be an
  injective resolution. The Comparison Theorem (see e.g.\/ \cite{Weibel},
  2.3.7) implies that $\mathcal{O}_{Y_{U}} \rightarrow
  \mathcal{J}_{U}\cx$ factors through $f_{U}^{\ast}\mathcal{O}_{U}
  \rightarrow   f_{u}^{\ast}\mathcal{I}_{U}\cx$, so that we obtain two
  quasi-isomorphisms
  $$\mathcal{O}_{Y_{U}} \rightarrow f_{U}^{\ast}\mathcal{I}_{U}\cx
  \rightarrow \mathcal{J}_{U}\cx.$$
  By the construction of derived functors, there are natural
  isomorphisms $\mathbb{A}_{U}(\mathcal{E}\cx) \cong \mathcal{P}_{U}^{\vee}
  \otimes f_{U}^{\ast} \mathcal{H}om_{U}(\mathcal{E}\cx,
  \mathcal{I}_{U}\cx)$ and
  $\mathbb{B}_{U}(\mathcal{E}\cx) \cong
  \mathcal{H}om_{Y_{U}}(\mathcal{P}_{U}\otimes f_{U}^{\ast}\mathcal{E}\cx,
  \mathcal{J}_{U}\cx)$, such that the choice of a different resolution leads
  to the isomorphism which is obtained by composing the above with the map
  which is induced by the unique homotopy equivalence between the
  resolutions.  
  We define $\eta_{U}(\mathcal{E}\cx)$ to be the composition of these
  isomorphisms with the following morphisms: 
  $$\begin{CD}
    \mathcal{P}_{U}^{\vee} \otimes f_{U}^{\ast}
    \mathcal{H}om_{U}(\mathcal{E}\cx, \mathcal{I}_{U}\cx)\\
    @VV{\eta_{1}}V\\
    \mathcal{P}_{U}^{\vee} \otimes
    \mathcal{H}om_{Y_{U}}(f_{U}^{\ast}\mathcal{E}\cx,
    f_{U}^{\ast}\mathcal{I}_{U}\cx) \\
    @VV{\eta_{2}}V\\
    \mathcal{P}_{U}^{\vee} \otimes
    \mathcal{H}om_{Y_{U}}(f_{U}^{\ast}\mathcal{E}\cx, \mathcal{J}_{U}\cx) \\
    @VV{\eta_{3}}V\\
    \mathcal{H}om_{Y_{U}}(\mathcal{P}_{U}\otimes f_{U}^{\ast}\mathcal{E}\cx,
    \mathcal{J}_{U}\cx).
  \end{CD}$$
  The morphism $\eta_{2}$ is induced by the quasi-isomorphism
  $f_{U}^{\ast}\mathcal{I}_{U}\cx \rightarrow
  \mathcal{J}_{U}\cx$. The morphisms $\eta_{1}$ and $\eta_{3}$ are the
  canonical ones, see for example \cite{EGA}, 0$_{\text{I}}$ 4.4.6 and 5.4.2.
  All these maps are natural with respect to $\mathcal{E}\cx\in\Dmcoh(U)$.
  Although, $\eta_{1}, \eta_{2}$ and $\eta_{3}$ are homomorphisms of complexes
  of $\mathcal{O}_{Y_{U}}$-modules, we are interested in
  $\eta_{U}(\mathcal{E}\cx)$ as a morphism in $\Dpcoh(Y_{U})$ only. As such,
  it does not depend on the choice of the resolutions $\mathcal{I}_{U}\cx$ and
  $\mathcal{J}_{U}\cx$. This follows from the Comparison Theorem and
  the fact that homotopy equivalent morphisms coincide as morphisms in the
  derived category.

  An important tool in our proof will be the following compatibility. If
  $U\subset T$ is an open subset, there exist isomorphisms of functors
  $\Dmcoh(T) \rightarrow \Dpcoh(Y_{U})$
  $$\mathbb{A}_{U}\circ i_{U}^{\ast} \cong j_{U}^{\ast} \circ \mathbb{A}_{T}
  \quad \text{ and } \quad
  \mathbb{B}_{U}\circ i_{U}^{\ast} \cong j_{U}^{\ast} \circ \mathbb{B}_{T}$$
  such that the diagram
  \begin{equation}\label{diag:restrict}
    \begin{CD}
      \mathbb{A}_{U}(i_{U}^{\ast}(\mathcal{E}\cx))
      @>{\eta_{U}(i_{U}^{\ast}(\mathcal{E}\cx))}>>
      \mathbb{B}_{U}(i_{U}^{\ast}(\mathcal{E}\cx))\\
      @VV{\cong}V @VV{\cong}V \\
      j_{U}^{\ast}\mathbb{A}_{T}(\mathcal{E}\cx)
      @>{j_{U}^{\ast}(\eta_{T}(\mathcal{E}\cx))}>>
      j_{U}^{\ast}\mathbb{B}_{T}(\mathcal{E}\cx)
    \end{CD}
  \end{equation}
  is commutative. This follows from our definition, because $j_{U}^{\ast}
  \circ f_{U}^{\ast} = f_{U}^{\ast} \circ i_{U}^{\ast}$, as well as
  $i_{U}^{\ast}\mathcal{H}om_{T}(\mathcal{F}\cx, \mathcal{G}\cx) =
  \mathcal{H}om_{U}(i_{U}^{\ast}\mathcal{F}\cx, i_{U}^{\ast}\mathcal{G}\cx)$
  and because the definition of $\eta$ does not depend on the choice of
  resolutions, as observed above.

  Next, we prove for any open $U\subset T$ that $\eta_{U}(\mathcal{O}_{U})$ is
  an isomorphism in $\Dpcoh(Y_{U})$. Because
  $\mathcal{E}xt_{U}^{i}(\mathcal{O}_{U}, \mathcal{O}_{U})=0$ for all $i>0$,
  the complex $\boldsymbol{R}\mathcal{H}om_{U}(\mathcal{O}_{U},
  \mathcal{O}_{U})$ is isomorphic to the sheaf $\mathcal{O}_{U}$ concentrated
  in degree zero. Hence, $\mathbb{A}_{U}(\mathcal{O}_{U}) \cong
  \mathcal{P}_{U}^{\vee}$. On the other hand,
  $\mathcal{E}xt_{Y}^{i}(\mathcal{P}, \mathcal{O}_{Y}) = 0$ for all $i>0$
  implies $\mathbb{B}_{U}(\mathcal{O}_{U}) \cong
  \mathcal{H}om(\mathcal{P}_{U}, \mathcal{O}_{U}) =
  \mathcal{P}_{U}^{\vee}$. By definition, $\eta_{U}(\mathcal{O}_{U})$ is the
  composition of the top row in the commutative diagram
  $$\begin{CD}
    \mathcal{P}_{U}^{\vee} \otimes f_{U}^{\ast}\mathcal{I}_{U}\cx
    @>>{\eta_{2}\circ \eta_{1}}>
    \mathcal{P}_{U}^{\vee} \otimes \mathcal{J}_{U}\cx
    @>>{\eta_{3}}>
    \mathcal{H}om_{Y_{U}}(\mathcal{P}_{U}, \mathcal{J}_{U}\cx)\\
    @AAA @AAA @AAA\\
    \mathcal{P}_{U}^{\vee} @= \mathcal{P}_{U}^{\vee} @= 
    \mathcal{H}om_{Y_{U}}(\mathcal{P}_{U}, \mathcal{O}_{Y_{U}}).
  \end{CD}$$
  Because the functor $\mathcal{P}_{U}^{\vee}\otimes f_{U}^{\ast}(\ARG)$ is
  exact, the first vertical arrow is a quasi-isomorphism. The third vertical
  arrow is a quasi-isomorphism as well, because
  $\mathcal{E}xt^{i}_{Y}(\mathcal{P}, \mathcal{O}_{Y}) = 0$ for all
  $i>0$. Hence, $\eta_{U}(\mathcal{O}_{U})$ is a quasi-isomorphism. Because
  $\mathbb{A}_{U}, \mathbb{B}_{U}$ and $\eta_{U}$ are compatible with finite
  direct sums, $\eta_{U}(\mathcal{O}_{U}^{\oplus k})$ is an isomorphism in
  $\Dpcoh(Y_{U})$ for all $k\ge 1$.

  If $\mathcal{E}$ is now an arbitrary locally free
  sheaf of finite rank on $T$, using diagram (\ref{diag:restrict}) we obtain
  that there exists an open covering of $T$ such that
  $j_{U}^{\ast}(\eta_{T}(\mathcal{E}))$ is an isomorphism for any open set
  $U\subset T$ in this covering. This implies that $\eta_{T}(\mathcal{E})$ is
  an isomorphism. The same statement is true for $\eta_{U}(\mathcal{E})$ if
  $\mathcal{E}$ is a locally free coherent sheaf on an arbitrary open subset
  $U\subset T$. 
  Because on any quasi-projective open subset $U\subset T$, any coherent sheaf
  has a resolution by coherent locally free sheaves, \cite{RD}, Prop.\/ I.7.1
  (iv), implies that $\eta_{U}(\mathcal{E})$ is an isomorphism if
  $\mathcal{E}$ is an arbitrary coherent sheaf on a quasi-projective open
  subset $U\subset T$. 

  We assumed $T$ to be of finite type over $\boldsymbol{k}$. This implies that
  there exists an open affine cover of $T$ by quasi-projective schemes. As
  before, using the diagram (\ref{diag:restrict}), we obtain that
  $\eta_{T}(\mathcal{E})$ is an isomorphism for any coherent sheaf
  $\mathcal{E}$ on $T$. The claim follows now from the lemma on way-out
  functors \cite{RD}, I.7.1. 
\end{proof}

Grothendieck-Verdier duality can be formulated in the following way (see
\cite{RD, Conrad}): if $f:X\rightarrow Y$ is a 
proper morphism between schemes which are of finite type over
$\boldsymbol{k}$, there exists a functor 
$$f^{!}:\Dpcoh(Y)\rightarrow \Dpcoh(X),$$ such that there is an isomorphism   
$$\boldsymbol{R}f_{\ast}\boldsymbol{R}\mathcal{H}om_{X}(\mathcal{F}\cx,
f^{!}\mathcal{G}\cx) \cong
\boldsymbol{R}\mathcal{H}om_{Y}
(\boldsymbol{R}f_{\ast}\mathcal{F}\cx,\mathcal{G}\cx),
$$
which is functorial in both arguments $\mathcal{F}\cx\in \Dmqc(X)$ and
$\mathcal{G}\cx \in \Dpcoh(Y)$. 
On the other hand, an object $\mathcal{I}\cx\in \Dbcoh(Y)$ of finite
injective dimension is called a dualising complex, if the functor
$\mathbb{D}_{Y}:= 
\boldsymbol{R}\mathcal{H}om_{Y}(\ARG, \mathcal{I}\cx)$ 
satisfies $\mathbb{D}_{Y}\circ\mathbb{D}_{Y}\cong\boldsymbol{1}$ on the
category $\Dcoh(Y)$.
Such a dualising complex exists on any scheme $Y$ which is of finite type over
$\boldsymbol{k}$.  
An important feature of the functor $f^{!}$ is that it respects the property
of a complex to be dualising. Hence, the functor $\mathbb{D}_{X}:= 
\boldsymbol{R}\mathcal{H}om_{X}(\ARG, f^{!}\mathcal{I}\cx)$ is
dualising and the duality theorem implies 
$$\boldsymbol{R}f_{\ast}\circ\mathbb{D}_{X} \cong
\mathbb{D}_{Y}\circ\boldsymbol{R}f_{\ast}.$$

If $\mathcal{I}\cx$ is a dualising complex on $Y$, for any integer $n$
and any invertible sheaf $\mathcal{L}$ on $Y$, the complex
$\mathcal{I}\cx[n] \otimes \mathcal{L}$ is a dualising complex as
well. Up to such changes, dualising complexes are unique. The notation
$\mathbb{D}_{Y}$ must be used with care, because such a functor depends on the
choice of a dualising complex on $Y$, but only up to a shift and a twist by an
invertible sheaf.

In order to apply the duality theorem, it is important to be able to calculate
$f^{!}$. This is particularly easy in the special case of a Cohen-Macaulay
morphism. If $Y$ is of finite type over $\boldsymbol{k}$ and $f:X\rightarrow
Y$ is a projective Cohen-Macaulay morphism whose fibres are 
of pure dimension $n$, we have for any $\mathcal{F}\cx\in \Dpcoh(Y)$ by
\cite{Conrad}, Theorem 4.3.2: 
$$f^{!}\mathcal{F}\cx\cong \omega_{f}[n]\dtens f^{\ast}\mathcal{F}\cx.$$ 
Here we denoted by $\omega_{f}$ the relative dualising sheaf, which is a
coherent $Y$-flat sheaf on $X$ and which satisfies 
$\omega_{f}[n]\cong f^{!}\mathcal{O}_{Y}$, see \cite{Conrad}, Theorem
3.5.1. This sheaf is locally free if and only if the morphism $f$ is
Gorenstein. The sheaf $\omega_{f}$ is compatible with base-change
(\cite{Conrad}, Theorem 3.6.1). This 
implies that its restriction to fibres is the dualising sheaf on the
fibre. In particular, if $f$ is a proper Gorenstein morphism all fibres of
which have trivial dualising sheaf and if $Y$ is reduced, then  
$$\omega_{f} \cong f^{\ast}\mathcal{A}$$
with an invertible sheaf $\mathcal{A}$ on $Y$.

In the formulation of the following theorem we use notation from Section
\ref{sec:fmt}.  
If $q:X\rightarrow S$ satisfies the assumptions ($\star$), there exists
an invertible sheaf $\mathcal{A}$ on $S$ such that $\omega_{q}\cong
q^{\ast}\mathcal{A}$. 
The sheaf $\mathcal{P}$ on $X\times_{S}X$ was defined in (\ref{eq:poincare})
and $\mathcal{M}$ is the sheaf introduced in Proposition
\ref{prop:i}, which satisfies $(\boldsymbol{1}_{X}\times i)^{\ast} \mathcal{P}
\cong \mathcal{P}^{\vee} \otimes \pi^{\ast}\mathcal{M}$.
Finally, $\mathcal{L}$ denotes an invertible sheaf on
$S$. Note that the functor $\boldsymbol{R}\mathcal{H}om_{X}(\ARG,
q^{\ast}\mathcal{L})$ is a dualising functor on $X$ only if $X$ is
Gorenstein. This case is considered in Corollary \ref{cor:alldual}.

\begin{theorem}\label{thm:alldual}
  Suppose $S, X$ and $q:X\rightarrow S$ satisfy condition ($\star$) from
  Section \ref{sec:fmt}. 
  Then, there exists an isomorphism in $\Dpcoh(X)$, which is functorial in
  $\mathcal{E}\cx\in\Dmcoh(X)$: 
  $$\boldsymbol{R}\mathcal{H}om_{X}(\FM_{\mathcal{P}}^{-}(\mathcal{E}\cx),
  q^{\ast}\mathcal{L}) 
  \cong
  i^{\ast}\FM_{\mathcal{P}}^{+}(
  \boldsymbol{R}\mathcal{H}om_{X}(\mathcal{E}\cx,q^{\ast}\mathcal{L})) 
  \otimes q^{\ast}(\mathcal{M}^{\vee}\otimes \mathcal{A})[1]$$  
\end{theorem}

\begin{proof}
  By definition, we have
  $$\boldsymbol{R}\mathcal{H}om_{X}(\FM_{\mathcal{P}}^{-}(\mathcal{E}\cx),
  q^{\ast}\mathcal{L}) 
  \cong
  \boldsymbol{R}\mathcal{H}om_{X}(\boldsymbol{R}\pi_{2\ast}(
  \mathcal{P}\otimes\pi_{1}^{\ast}\mathcal{E}\cx),q^{\ast}\mathcal{L})$$
  and apply Grothendieck-Verdier duality to obtain
  $$\cong
  \boldsymbol{R}\pi_{2\ast}\boldsymbol{R}\mathcal{H}om_{X\times_{S}X}(
  \mathcal{P}\otimes\pi_{1}^{\ast}\mathcal{E}\cx,
  \pi_{2}^{!}q^{\ast}\mathcal{L}).$$
  Using \cite{Conrad}, Theorem 3.6.1, we see
  $$\pi_{2}^{!}q^{\ast}\mathcal{L} \cong \omega_{\pi_{2}}[1] \otimes
  \pi_{2}^{\ast}q^{\ast}\mathcal{L} \cong \pi_{1}^{\ast}\omega_{q}[1] \otimes
  \pi_{2}^{\ast}q^{\ast}\mathcal{L} \cong
  \pi^{\ast}(\mathcal{A}\otimes \mathcal{L})[1].$$ 
  Therefore, the above is isomorphic to  
  $$\boldsymbol{R}\pi_{2\ast}\big(\boldsymbol{R} \mathcal{H}om_{X\times_{S}X}
  (\mathcal{P}\otimes\pi_{1}^{\ast}\mathcal{E}\cx, \mathcal{O}_{X\times_{S}X})
  \otimes \pi^{\ast}(\mathcal{A}\otimes \mathcal{L})[1]\big).$$ 
  Applying Lemma \ref{lem:dual} and using 
  $\pi^{\ast}\mathcal{L} \cong \pi_{1}^{\ast}q^{\ast}\mathcal{L}$,
  we see that this is functorially isomorphic to
  $$\boldsymbol{R}\pi_{2\ast}\big(\mathcal{P}^{\vee} \otimes \pi_{1}^{\ast}
  \boldsymbol{R} \mathcal{H}om_{X}(\mathcal{E}\cx, q^{\ast}\mathcal{L}) \otimes
  \pi^{\ast}\mathcal{A}\big)[1].$$ 
  Now we use $\pi^{\ast}\mathcal{A} \cong \pi_{2}^{\ast}q^{\ast}\mathcal{A}$, 
  $\pi_{1}\circ(\boldsymbol{1}_{X}\times i) = \pi_{1}$ and
  $\mathcal{P}^{\vee} \cong (\boldsymbol{1}_{X}\times i)^{\ast}\mathcal{P}
  \otimes \pi^{\ast}\mathcal{M}^{\vee}$ with $\mathcal{M}$ being the
  invertible sheaf on $S$ introduced in Proposition \ref{prop:i}. This allows
  us to write the above expression as
  $$\boldsymbol{R}\pi_{2\ast}\big(
  (\boldsymbol{1}_{X}\times i)^{\ast}\mathcal{P} \otimes
  (\boldsymbol{1}_{X}\times i)^{\ast} \pi_{1}^{\ast}
  \boldsymbol{R} \mathcal{H}om_{X}(\mathcal{E}\cx, q^{\ast}\mathcal{L}) \otimes
  \pi_{2}^{\ast}q^{\ast}(\mathcal{M}^{\vee}\otimes\mathcal{A})\big)[1].$$
  Finally, we apply the projection formula, followed by 
  $\boldsymbol{R}\pi_{2\ast} \circ (\boldsymbol{1}_{X}\times i)^{\ast}
  \cong  i^{\ast} \circ \boldsymbol{R}\pi_{2\ast}$, which is base-change with
  the isomorphism $i$, and obtain the desired isomorphism
  $$\boldsymbol{R}\pi_{2\ast}\big(
  (\boldsymbol{1}_{X}\times i)^{\ast}(\mathcal{P}
  \otimes \pi_{1}^{\ast}
  \boldsymbol{R} \mathcal{H}om_{X}(\mathcal{E}\cx, q^{\ast}\mathcal{L}))\big)
  \otimes 
  q^{\ast}(\mathcal{M}^{\vee}\otimes\mathcal{A})[1]$$ 
  $$\cong i^{\ast}\FM_{\mathcal{P}}^{+}(
  \boldsymbol{R}\mathcal{H}om_{X}(\mathcal{E}\cx, q^{\ast}\mathcal{L})) \otimes
  q^{\ast}(\mathcal{M}^{\vee}\otimes \mathcal{A})[1].$$
\end{proof}

\begin{corollary}\label{cor:alldual}
  In addition to the notation and assumptions of Theorem \ref{thm:alldual},
  suppose $S$ is Gorenstein and connected. Then, there exist isomorphism
  \begin{align*}
    \mathbb{D}_{X} \circ \FM_{\mathcal{P}}^{-} &\cong
  [1]\circ\mathbb{T}_{\mathcal{M}^{\vee}\otimes\mathcal{A}} \circ
  i^{\ast}\circ\FM_{\mathcal{P}}^{+}\circ\mathbb{D}_{X} \quad\text{and}\\
  \mathbb{D}_{X} \circ \FM_{\mathcal{P}}^{\mathsf{b}} &\cong
  [1]\circ\mathbb{T}_{\mathcal{M}^{\vee}\otimes\mathcal{A}} \circ
  i^{\ast}\circ\FM_{\mathcal{P}}^{\mathsf{b}}\circ\mathbb{D}_{X},
  \end{align*}
  where $\mathbb{T}_{\mathcal{M}^{\vee}\otimes\mathcal{A}}$ denotes the tensor
  product functor with the locally free sheaf
  $q^{\ast}(\mathcal{M}^{\vee}\otimes\mathcal{A})$ and $\mathbb{D}_{X}$ is a
  dualising functor on $X$ of the form
  $\boldsymbol{R}\mathcal{H}om_{X}(-,q^{\ast}\mathcal{L})$, with an invertible 
  sheaf $\mathcal{L}$ on $S$.
\end{corollary}

\begin{proof}
  Because the scheme $S$ and the morphism $q$ are Gorenstein, the scheme $X$
  is Gorenstein as well. Hence, any shift of a locally free sheaf on $X$ is a
  dualising complex. In particular, for any locally free sheaf $\mathcal{L}$
  on $S$, the functor $\boldsymbol{R}\mathcal{H}om_{X}(-,q^{\ast}\mathcal{L})$
  is dualising on $X$. The claim is just a reformulation of the theorem.
  The statement involving $\FM_{\mathcal{P}}^{\mathsf{b}}$ is obtained by
  restricting the other one, because $X$ is Gorenstein, hence $\mathbb{D}_{X}:
  \Dbcoh(X) \rightarrow \Dbcoh(X)$.
\end{proof}

In the special case $S = \Spec(\boldsymbol{k})$, the invertible sheaves
$\mathcal{M}$ and $\mathcal{A}$ are automatically trivial, so that
$\mathbb{T}_{\mathcal{M}^{\vee}\otimes\mathcal{A}} \cong
\boldsymbol{1}$. Because any  
irreducible projective curve of arithmetic genus one is automatically
Gorenstein with trivial canonical sheaf, we obtain the following corollary,
which generalises Mukai's result \cite{Mukai}, (3.8). We denote
$\mathbb{D}_{\boldsymbol{E}} =
\boldsymbol{R}\mathcal{H}om_{\boldsymbol{E}}(-,\mathcal{O}_{\boldsymbol{E}})$.

\begin{corollary}\label{cor:dual}
  If $\boldsymbol{E}$ is an irreducible projective curve of arithmetic genus
  one, then: 
  \begin{align*}
    \mathbb{D}_{\boldsymbol{E}} \circ \FM_{\mathcal{P}}^{-} &\cong
    [1]\circ i^{\ast}\circ\FM_{\mathcal{P}}^{+}\circ\mathbb{D}_{\boldsymbol{E}}
    \quad\text{and}\\
    \mathbb{D}_{\boldsymbol{E}} \circ \FM_{\mathcal{P}}^{\mathsf{b}} &\cong
   [1]\circ i^{\ast}\circ \FM_{\mathcal{P}}^{\mathsf{b}} \circ
    \mathbb{D}_{\boldsymbol{E}}. 
  \end{align*}
\end{corollary}

This is a generalisation of \cite{BurbanKreussler}, Theorem 6.11. 
To see this, recall that we studied in \cite{BurbanKreussler} the
functor $\mathbb{F}:= [1] \circ \FM_{\mathcal{P}}^{\mathsf{b}}$.
Moreover, for any coherent torsion sheaf $\mathcal{F}$ we called
$\mathbb{M}(\mathcal{F}) = \mathcal{E}xt^{1}(\mathcal{F},\mathcal{O})$ the
Matlis dual of $\mathcal{F}$. 
Note that we have $\mathcal{E}xt^{i}(\mathcal{F},\mathcal{O})=0$ for such
$\mathcal{F}$ and $i\ne 1$. This implies
$\mathbb{D}_{\boldsymbol{E}}(\mathcal{F}) \cong \mathbb{M}(\mathcal{F})[-1]$.
For any semi-stable torsion free sheaf $\mathcal{E}$ of degree zero on
$\boldsymbol{E}$ we have shown in \cite{BurbanKreussler}, Theorem 6.11, that
there is an 
isomorphism $\mathbb{M}(\mathbb{F}(\mathcal{E})) \cong
i^{\ast}\mathbb{F}(\mathcal{E}^{\vee})$. 
Now, observe that we have $\mathcal{E}xt^{i}(\mathcal{E},\mathcal{O})=0$ for
any torsion free sheaf $\mathcal{E}$ on $\boldsymbol{E}$ and all $i>0$. 
Hence, $\mathbb{D}_{\boldsymbol{E}}(\mathcal{E}) \cong \mathcal{E}^{\vee}$. 
Corollary \ref{cor:dual} together with $\mathbb{D}_{\boldsymbol{E}}\circ[-1]
\cong [1]\circ\mathbb{D}_{\boldsymbol{E}}$ gives, therefore, 
$$\mathbb{M}(\mathbb{F}(\mathcal{E})) \cong
\mathbb{D}_{\boldsymbol{E}}(\mathbb{F}(\mathcal{E}))[1] \cong
\mathbb{D}_{\boldsymbol{E}}(\FM_{\mathcal{P}}^{\mathsf{b}}(\mathcal{E})) \cong
i^{\ast}\FM_{\mathcal{P}}^{\mathsf{b}}
(\mathbb{D}_{\boldsymbol{E}}(\mathcal{E}))[1]$$ 
$$\cong i^{\ast} \mathbb{F}(\mathbb{D}_{\boldsymbol{E}}(\mathcal{E}))\cong
i^{\ast}\mathbb{F}(\mathcal{E}^{\vee}).$$

\section{Applications}
\label{sec:appl}

The purpose of this section is to show through simple examples how the
results of the previous sections can be used to explicitly construct
interesting families of sheaves.
In order to be able to study sheaves on the fibres, we need the following
compatibility of Fourier-Mukai transforms with restriction functors.

\begin{lemma}\label{lem:BW}
  Suppose $q: X \rightarrow S$ satisfies condition ($\star$) from Section
  \ref{sec:fmt}, let $s: \Spec(\boldsymbol{k}) \rightarrow S$ be an arbitrary 
  point and denote by $i_{s}:X_{s}\rightarrow X$ the embedding of the fibre of
  $q$ over $s$. Let $\mathcal{P}$ be the sheaf (\ref{eq:poincare}) on
  $X\times_S X$. By $\mathcal{P}_s$ we denote its restriction to $X_s\times
  X_s$. If $\mathcal{E}$ is an $S$-flat coherent sheaf on $X$, there exists an
  isomorphism
  \begin{equation}
    \label{eq:BW}
    \boldsymbol{L}i_s^\ast \FM_{\mathcal{P}}^{\mathsf{b}}(\mathcal{E})
    \cong \FM_{\mathcal{P}_s}^{\mathsf{b}} i_s^\ast(\mathcal{E}).
  \end{equation}
  Moreover, if $\FM_{\mathcal{P}_s}^{\mathsf{b}} i_s^\ast(\mathcal{E})$ is a
  sheaf for all $s\in S$, $\FM_{\mathcal{P}}^{\mathsf{b}}(\mathcal{E})$ is an
  $S$-flat sheaf on $X$ and we can replace $\boldsymbol{L}i_s^\ast$ by
  $i_s^\ast$ in this formula.
\end{lemma}

\begin{proof}
  Because $\mathcal{E}$ and $\mathcal{P}$ are both $S$-flat, the sheaf
  $\mathcal{P}\otimes \pi_{1}^{\ast}\mathcal{E}$ is $\pi_{2}$-flat as
  well. Hence, we have a base-change isomorphism (see \cite{Mumford}, II \S5)
  $$\boldsymbol{L}i_{s}^{\ast}\boldsymbol{R}\pi_{2\ast}
  (\mathcal{P}\otimes \pi_{1}^{\ast}\mathcal{E})
  \cong
  \boldsymbol{R}p_{2\ast}(i_{s}\times i_{s})^{\ast}
  (\mathcal{P}\otimes \pi_{1}^{\ast}\mathcal{E}),$$
  which implies (\ref{eq:BW}). The final statement follows from
  \cite{BridgelandEquiv}, Lemma 4.3. 
\end{proof}

We  conclude this article by presenting two examples of fibrewise 
semi-stable torsion free sheaves on an elliptic fibration 
$X \rightarrow S$.
The main idea of the construction is to apply the relative Fourier-Mukai
transform $\FM_{\mathcal{P}}^{\mathsf{b}}$ to an $S$-flat family of torsion
sheaves on $X$. The result is a flat family of semi-stable torsion free
sheaves. This is a special case of the spectral cover construction \cite{FMW,
DonagiPantev}. 
Moreover, using  our previous results  about Fourier-Mukai transforms on 
Weierstra\ss{}  cubics \cite{BurbanKreussler}, we are able to describe
explicitly the restriction onto the singular fiber in terns of \'etale
coverings of nodal cubic curves. 

Let $X\subset \mathbb{P}^{2}\times \mathbb{A}^{1}$ be the family of
irreducible cubics, defined by the polynomial
$$F_{t}(x,y,z) = y^{2}z - x^{3} - x^{2}z - t(1-t)xz^{2} + t^{2}z^{3},$$
where $(x:y:z)$ are homogeneous coordinates on $\mathbb{P}^{2}$ and $t$ is an
affine coordinate on $\mathbb{A}^{1}$. The discriminant of this family is a
polynomial in $t$ of degree six which vanishes at $t=0$.
The general fibre of the projection $X\rightarrow\mathbb{A}^{1}$ is
smooth. The fibre $X_{0}$ over $t=0$ is given by the polynomial $y^{2}z -
x^{3} - x^{2}z$, hence is an irreducible cubic with a node at $(x:y:z) =
(0:0:1)$. 
This family has a section $\sigma:\mathbb{A}^{1} \rightarrow X$, which is
given by $\sigma(t)=((0:1:0),t)\in X$.

\begin{example}\label{expl:familyone}
  For any $\lambda\in\boldsymbol{k}^{\times}$, we denote by $C_{\lambda}$ the
  intersection of $X$ with the 
  surface which is given in $\mathbb{P}^{2}\times \mathbb{A}^{1}$ by the
  equation $$(1+\lambda)y - (1-\lambda)x=0.$$
  We first study the case $\lambda\ne-1$. In this case, $C_{\lambda}$ is
  disjoint to the set where $z=0$. Therefore, it is sufficient to work on the
  affine set $z\ne0$, where the curve $C_{\lambda}$ is defined by the
  ideal
  $$I_{\lambda} = \langle (1+\lambda)y - (1-\lambda)x,\;
  y^{2} - x^{3} - x^{2} - t(1-t)x + t^{2}\rangle.$$
  If $t\in\boldsymbol{k}$ is arbitrary, the module
  $\boldsymbol{k}[x,y]/I_{\lambda}$ is of length three. For generic $t$ 
  it consists of three simple points.
  Note that the sheaf $\mathcal{O}_{C_{\lambda}}$ is 
  flat over $\mathbb{A}^{1}$, because $t$ is not a zero-divisor of
  $\boldsymbol{k}[x,y,t]/I_{\lambda}$.

  The fibre of $\mathcal{O}_{C_{\lambda}}$ over $t=0$ is given by the ideal
  \begin{align*}
    &\langle (1+\lambda)y - (1-\lambda)x,\; y^{2} - x^{3} - x^{2}\rangle\\
    =&\left\langle (1+\lambda)y - (1-\lambda)x,\;
    x^{2}\left(\frac{4\lambda}{(1+\lambda)^{2}} +x\right)\right\rangle
    \subset \boldsymbol{k}[x,y].
  \end{align*}
  Thus, the support of $i_{0}^{\ast}\mathcal{O}_{C_{\lambda}}$ consists of the
  two points $(0:0:1)$ and $(4\lambda(1+\lambda) : 4\lambda(1-\lambda) :
  -(1+\lambda)^{3})$ in $\mathbb{P}^{2}$.
  At the second of them, $i_{0}^{\ast}\mathcal{O}_{C_{\lambda}}$ is of length
  one. In order to 
  understand $i_{0}^{\ast}\mathcal{O}_{C_{\lambda}}$ at the singular point, we
  look at the completion $\widehat{M}$ at the point $(x,y)=(0,0)$ of
  $$M=\boldsymbol{k}[x,y]/ \left\langle (1+\lambda)y - (1-\lambda)x,\; 
  x^{2}\left(4\lambda +(1+\lambda)^{2}x\right)\right\rangle.$$
  As we always assume $\lambda\ne0$, the factor behind $x$ is a unit in
  $\widehat{R} = \boldsymbol{k}[[x,y]]/ \langle y^{2}-x^{2}-x^{3}\rangle$,
  hence $\widehat{M} = \boldsymbol{k}[[x,y]]/ \left\langle (1+\lambda)y -
  (1-\lambda)x,\; x^{2}\right\rangle$.
  If we let $\widetilde{x}:=x\sqrt{1+x}\in\widehat{R}$ and define 
  $\xi=\widetilde{x}-y, \eta=\widetilde{x}+y$, we obtain an isomorphism
  $\widehat{R} \cong \boldsymbol{k}[[\xi,\eta]]/\langle \xi\eta\rangle$.
  
  Note that $\sqrt{1+x}= 1+x/2-x^{2}/8+\ldots$, as given by the
  binomial series, is a unit in $\widehat{R}$. Because
  $\widetilde{x}-x\in\langle x^{2}\rangle \subseteq \widehat{R}$ and $x^{2}=0$
  in $\widehat{M}$, we obtain in $\widehat{M}$
  $$\lambda\eta - \xi =
  \lambda(y+\widetilde{x}) +  y - \widetilde{x}
  = (1+\lambda)y - (1-\lambda)x = 0.$$ 
  This implies $\widehat{M} \cong \widehat{R}/\langle \lambda\eta -\xi\rangle$.
  Using the notation from \cite{BurbanKreussler}, \S4, this module is
  identified with the band module $\mathcal{M}((1,1),1,\lambda)$ in the
  classification of Gelfand and Ponomarev. It has length two, as expected.

  The functor
  $\mathbb{F}=[1]\circ\FM_{\mathcal{P}_{0}}^{\mathsf{b}}$ was shown in 
  \cite{BurbanKreussler} to satisfy $\mathbb{F}\circ\mathbb{F} \cong
  i^{\ast}[1]$ and $\mathbb{F}(\mathcal{B}((1,-1),1,\lambda)) \cong
  \mathcal{M}((1,1),1,\lambda)$. Here, $i$ is an involution on $X_{0}$ and
  $\mathcal{B}((1,-1),1,\lambda)$ a semi-stable
  indecomposable vector bundle of rank two and degree zero on $X_{0}$.
  This implies
  $$\FM_{\mathcal{P}_{0}}^{\mathsf{b}}(i_{0}^{\ast}\mathcal{O}_{C_{\lambda}})
  \cong 
  \mathcal{B}((1,-1),1,\lambda)\oplus \mathcal{L}_{\lambda},$$
  where $\mathcal{L}_{\lambda}$ is a line bundle on $X_{0}$ which depends on
  $\lambda$.

  As we have seen above, for generic $t\in\mathbb{A}^{1}$ the fibre of
  $C_{\lambda}$ over $t$ consists of three points of length one. This implies
  that
  $\FM_{\mathcal{P}_{t}}^{\mathsf{b}}(i_{t}^{\ast}\mathcal{O}_{C_{\lambda}})$
  is the direct sum of three line bundles. 
  Using Lemma \ref{lem:BW}, we obtain that
  $\FM_{\mathcal{P}}^{\mathsf{b}}(\mathcal{O}_{C_{\lambda}})$ is a coherent
  $\mathbb{A}^{1}$-flat sheaf on $X$, which is the direct sum of three line
  bundles on the generic fibre, but which has a direct summand isomorphic to
  the indecomposable vector bundle $\mathcal{B}((1,-1),1,\lambda)$ if
  restricted to the singular fibre $X_{0}$. 

  In the case $\lambda=-1$, the situation is very similar. The main difference
  is that the image of the section $\sigma$ is a component of $C_{-1}$. On the
  affine open set $z\ne0$, we obtain $I_{-1} = \langle x,\; y^{2} +
  t^{2}\rangle$ and, for any fixed $t\in\boldsymbol{k}$, the module
  $\boldsymbol{k}[x,y]/I_{-1}$ has length two. This module  
  consists of two simple points for all $t\ne0$.
  Hence, the component of $C_{-1}$, which is not supported at the
  image of the section $\sigma$, is of degree two over $\mathbb{A}^{1}$. 
  The restriction of this component to the fibre $X_{0}$ is supported at the
  singular point. Its completion is isomorphic to $\widehat{R}/\langle
  \xi+\eta\rangle$. 
  This is the band module $\mathcal{M}((1,1),1,-1)$ about which we know
  $\mathbb{F}(\mathcal{B}((1,-1),1,-1)) \cong \mathcal{M}((1,1),1,-1)$.
  
  Using Lemma \ref{lem:BW}, this implies that
  $\FM_{\mathcal{P}}^{\mathsf{b}}(\mathcal{O}_{C_{-1}})$ splits 
  into a direct sum of two $\mathbb{A}^{1}$-flat families
  $\mathcal{A}_{1}\oplus\mathcal{A}_{2}$. The restriction of 
  $\mathcal{A}_{1}$ on each fibre is trivial.
  However, $\mathcal{A}_{2}$ is a coherent
  $\mathbb{A}^{1}$-flat sheaf on $X$, which is the direct sum of two line
  bundles on the generic fibre, but which is isomorphic to the indecomposable
  vector bundle $\mathcal{B}((1,-1),1,-1)$ if restricted to the singular fibre
  $X_{0}$.  
\end{example}

\begin{example}\label{expl:familytwo}
  In this example we benefit from the special choice of our family of
  cubics, which can be written in the form
  $$F_{t}(x,y,z)=(y^{2}-x^{2})z-(x-tz)(x^{2}+txz+tz^{2}).$$
  We define $C\subset\mathbb{P}^{2}\times\mathbb{A}^{1}$ by the two equations
  $$x+y=0\quad  \text{ and }\quad x^{2}+txz+tz^{2}=0.$$
  Clearly, $C\subset X$ and the support of $C$ 
  is disjoint to the set given by $z=0$. Therefore, we restrict our attention
  to the affine open set $z\ne0$, in which $C$ is given by the ideal 
  $$I=\langle x+y,\; x^{2}+tx+t\rangle\subset\boldsymbol{k}[x,y,t].$$
  For any $t\in\boldsymbol{k}$, the vector space $\boldsymbol{k}[x,y]/I$ is of
  dimension two. Moreover, $\boldsymbol{k}[x,y,t]/I$ is a flat
  $\boldsymbol{k}[t]$-module, because $t$ is not a zero-divisor of it. Hence,
  the sheaf $\mathcal{O}_{C}$ is a flat family of torsion sheaves of length
  two. For generic $t$, this module is supported at two distinct
  points whose coordinates $(x,y)=(x,-x)$ satisfy $x^{2}+tx+t=0$.
  
  The fibre of $\mathcal{O}_{C}$ over $t=0$ is isomorphic to the module
  $$M=\boldsymbol{k}[x,y]/\langle x+y,\; x^{2}\rangle$$
  over the ring $R=\boldsymbol{k}[x,y]/\langle y^{2}-x^{2}-x^{3}\rangle$, 
  and is supported at $(x,y)=(0,0)$ only.
  Using $\widetilde{x}=x\sqrt{1+x}\in\widehat{R}$ and $\xi=\widetilde{x}-y, 
  \eta=\widetilde{x}+y$ as before, we obtain $\widetilde{x}=x$ in
  $\widehat{M}$, the completion of $M$ at $(0,0)$. This implies $\eta=0$ and
  $\xi^{2}=4\widetilde{x}^{2}=4x^{2}=0$ in $\widehat{M}$, so that we obtain
  $$\widehat{M} \cong \widehat{R}/\langle\xi^{2}, \eta\rangle.$$
  This module is called a string module in the classification of Gelfand,
  Ponomarev and was denoted $\mathcal{N}(0()1)$ in
  \cite{BurbanKreussler}. It has length two as expected.
  As in the previous example, we use the calculation from
  \cite{BurbanKreussler}, which shows 
  $\mathbb{F}(\mathcal{S}(0,-1)) \cong \mathcal{N}(0()1)$, where
  $\mathcal{S}(0,-1)$ is a semi-stable torsion free coherent sheaf of rank
  two and degree zero on $X_{0}$. Therefore, we obtain
  $$\FM_{\mathcal{P}_{0}}^{\mathsf{b}}(i_{0}^{\ast}\mathcal{O}_{C}) \cong
  \mathcal{S}(0,-1)$$
  and $\FM_{\mathcal{P}}^{\mathsf{b}}(\mathcal{O}_{C})$ is an
  $\mathbb{A}^{1}$-flat family of coherent sheaves on $X$ whose general fibre
  is the direct sum of two line bundles. The restriction to the singular fibre
  $X_{0}$, however, is indecomposable and not locally free.
\end{example}

\end{document}